\documentclass{siamart1116}

\usepackage{lipsum}
\usepackage{amsfonts}
\usepackage{graphicx}
\usepackage{epstopdf}
\usepackage{algorithmic}
\ifpdf
  \DeclareGraphicsExtensions{.eps,.pdf,.png,.jpg}
\else
  \DeclareGraphicsExtensions{.eps}
\fi

\numberwithin{theorem}{section}

\newcommand{\TheTitle}{Convergence characteristics of
  the generalized residual cutting method} 
\newcommand{\TheAuthors}{T. Abe, A. T.  Chronopoulos}

\headers{\TheTitle}{\TheAuthors}

\title{{\TheTitle}\thanks{Submitted to the editors DATE.}
}
\author{T. Abe\thanks{Advanced Science \& Intelligence Research Institute,
  Chiyoda-ku, Tokyo 101-0047, Japan
  }
    (\email{abe.toshihiko24@gmail.com})
  \and
  A. T.  Chronopoulos\thanks{
    Department of Computer Science, University of Texas,  San
    Antonio, TX 78249, USA
    and
    (Visiting Faculty) Department of Computer Science,
    University of  Patras, Greece (\email{antony.tc@gmail.com})
    }
}

\usepackage{amsopn}

\ifpdf
\hypersetup{
  pdftitle={\TheTitle},
  pdfauthor={\TheAuthors}
}
\fi

\newcommand{\eq}{\begin{equation}}
\newcommand{\qe}{\end{equation}}
\newcommand{\bgc}{\begin{center}}
\newcommand{\edc}{\end{center}}

\begin{document}

\maketitle

% REQUIRED
\begin{abstract}
  The residual cutting (RC) method has been proposed
  for efficiently solving linear equations obtained from
  elliptic partial differential equations.
  Based on the RC, we have introduced the generalized residual cutting (GRC)
  method, which can be applied to general sparse matrix problems.
  In this paper,
  we study the mathematics of the GRC algorithm and
  and prove it is a Krylov subspace method.
  Moreover, we show that it is deeply related to the conjugate residual (CR) method
  and that GRC becomes equivalent to CR for symmetric matrices.
  Also, in numerical experiments,
  GRC shows more robust convergence and needs less memory compared to GMRES,
  for significantly larger matrix sizes,

\end{abstract}



\section{Introduction}

The residual cutting (RC) method has been proposed
for efficiently solving linear equations obtained from
elliptic partial differential equations\cite{rc}.
It is an iterative method and for each iteration step,
its inner solver applies a relaxation method such as
successive over-relaxation (SOR)
to obtain an approximate solution.
Based on the principle of residual minimization, RC accelerates convergence of the inner solver.
However, its application is limited to linear problems
with diagonal-dominant matrices in general, for which convergence of a relaxation method
such as SOR is guaranteed.

We applied the RC method to coupled perturbed (CP) equations\cite{cprc},
to achieve robust convergence of iterative calculation.
Since relaxation methods are not applicable to CP equations in general due to reasons below,
some modification to RC was necessary.
Each iterative calculation in solving a CP equation is linear and
can be expressed by matrix-vector operations.
The matrix size is
proportional to the square of the number of orbitals
and becomes extremely large in general. Moreover, its elements are not directly accessible
in practice\cite{cprc}.
Since the matrix is neither sparse nor diagonally dominant,
iterative calculation sometimes diverges.
This is the main reason that we applied the RC method
to prevent divergence of the iterative solution.
However, direct application of the RC method was not possible
because it adopts a relaxation method in its inner solver, which needs
diagonal dominance and also direct access to matrix elements.
To overcome this problem, we replaced the inner solver in RC with
a matrix-vector operation\cite{cprc} (CPRC method).

The conjugate residual (CR) Krylov iterative  method for solving linear
systems has been studied  (or  presented) by several authors e.g.
\cite{ref1}-\cite{ref5}, \cite{ref14}, \cite{ref15}.
CR is a variant of the conjugate gradient Krylov iterative  method that
also applies to 
linear systems where the matrix of coefficients  is nonsymmetric and it has
positive definite symmetric part. 

We have introduced the generalized residual cutting (GRC) method \cite{grc}
as a Krylov subspace method. however, with no mathematical proof
based on matrix polynomials that it is a Krylov method.

In this paper, we study the mathematics of the GRC algorithm
and prove that it is a Krylov subspace method.
Also, we will show that it is deeply related to the conjugate residual (CR) method.
In fact, GRC becomes equivalent to CR for symmetric matrices.

\section{The residual cutting method}
\subsection{Original RC}

Here we cite a brief description of the RC method from our previous paper
\cite{cprc}, which the new method is based on.
We deal with the solution of a linear system of equations
\eq
\bf{H} \bf{U} = \bf{b}
\qe
are summarized as follows,
where $\bf{H}$ is the coefficient matrix, $\bf{U}$ is the solution vector
and $\bf{b}$ is the right hand side known vector.

%\newpage

Our analysis leads to the algorithm in \cref{alg:rcmethod}.

\begin{algorithm}
  \caption{RC method}
  \label{alg:rcmethod}
  \begin{algorithmic}
    \STATE{Set initial value for ${\bf U}^0$}
    \FOR{$m=0$ Until Convergence}
    \STATE{${\bf r}^m = {\bf b} - {\bf H} {\bf U}^m$}
    \STATE{Compute ${\bf \Psi}^m$ for ${\bf H} {\bf \Psi}^m = {\bf r}^m$}
    \STATE{Select $\alpha_k (k=1 \cdots L)$
      that minimize \\ \hspace{10mm}
      $||{\bf r}^{m+1}|| = ||{\bf b} - {\bf H} ({\bf U}^m + \phi^m)|| = ||r^m - {\bf H} \phi^m||$
      \\ \hspace {10mm} in $\phi^m = \alpha_1 {\bf \Psi}^m + \sum_{k=2}^{L} \alpha_k \phi^{m-k+1}$}
    \STATE{${\bf U}^{m+1} = {\bf U}^{m} + \phi^m$}
    \STATE{$r^{m+1} = r^m - {\bf H} \phi^m$}
    \ENDFOR
    \RETURN $\phi^m$
  \end{algorithmic}
\end{algorithm}

\subsection{The coupled perturbed RC (CPRC) method}
The generalized residual cutting (GRC) method has its origin in
our previous study of applying the RC method
to a coupled perturbed (CP)equation\cite{cprc}.
Each iterative calculation in solving a CP equation is linear and
can be expressed by matrix-vector operations.
The matrix size is
proportional to the square of the number of orbitals
and extremely large in general. Moreover, its elements are not directly accessible in practice.
Since the matrix is neither sparse nor diagonally dominant,
% it is impossible to keep the whole matrix
%in memory or even in a file, and a preconditioner is not applicable.
%Moreover, since it is in general not diagonally dominant,
iterative calculation sometimes diverges.
This is the main reason that we applied the RC method
to achieve robust convergence of the iterative solution.
However, direct application of the RC method was not possible
because it adopts a relaxation method in its inner solver, which needs
diagonal dominance and also direct access to matrix elements.
In order to apply the RC method,
we replaced the inner solver with a matrix-vector operation
\eq \label{eq.cprcpsi}
{\bf \Psi}^m = ({\bf 1 - H}) \phi^{m-1} + {\bf r}^m
\qe
Eq. (\ref{eq.cprcpsi}) is the matrix-vector expression of
the eq. (A.1) in \cite{cprc}, where it is called 'simple damping'.
The $\phi^{m-1}$ corresponds to the current unknowns ${\bf P}^{(1,0)}$
in \cite{cprc} and ${\bf \Psi}^m$ corresponds to the updated ${\bf P}^{(1,0)}$.
Convergence of the simple damping means ${\bf \Psi}^m = \phi^{m-1}$.
Then Eq. (\ref{eq.cprcpsi}) becomes ${\bf H} {\bf \Psi}^m = {\bf r}^m$,
which means that the equation for ${\bf \Psi}^m$ is solved.

Although the matrix $\bf{H}$ is not directly accessible, an equivalent
operation can be done with linear operations in a single iteration,
which consist of calculation of perturbed molecular orbital integrals.
%Note that elements of the matrix $H$ are not directly
%accessible since the whole matrix is too large to be kept in memory,
%or even in a file. For that reason, each elements are newly calculated
%for each iteration, as is the common case which is so-called 'direct method'.
The matrix-vector operation in (2.2) generates a new vector
and will be used for residual minimization.

This method of applying RC to CP equations will be denoted by 'CPRC' hereafter.
The difference between CPRC and the original RC is that in CPRC, the new vector does not
need to be an approximate solution and it will be used to
construct a new vector subspace.
On the other hand, in the original RC, it must be an approximate
solution used for accelerating convergence.

\subsection{The generalized residual cutting (GRC) method}

Calculation of \\ CPRC is identical to the Algorithm 2.1,
except for using (2.2) for 
computing ${\bf \Psi}^m$.
We will define the generalized residual cutting (GRC) method,
which includes CPRC as a special case.

In the conventional RC method, $\Psi^m$ is a temporary approximate solution
obtained by the inner solver, which usually adopts a relaxation method
such as successive over-relaxation (SOR).
When a relaxation method is used (here we show the case of Gauss-Seidel method
as an example), its single operation can be written as
\eq\label{eq.relax}
{\bf x}^i = ({\bf H}_{\rm L}+{\bf H}_{\rm D})^{-1}(-{\bf H}_{\rm U} {\bf x}^{i-1} + {\bf r}^m) \equiv f({\bf x}^{i-1})
\qe
where ${\bf H}_{\rm L}$, ${\bf H}_{\rm U}$ and ${\bf H}_{\rm D}$ are the lower-triangle, upper-triangle and
the diagonal components of the coefficient matrix.
In the RC method, the inner solver applies the above single operation for $n$ times,
then the temporary approximate solution becomes
\eq
{\bf \Psi}^m_{\scriptscriptstyle \rm RC} = f^n({\bf x}_0) = F({\bf H},{\bf r}^m)
\qe
where ${\bf x}_0$ is an initial vector (typically zero) for the relaxation method.
The last term means that ${\bf \Psi}^m_{\scriptscriptstyle \rm RC}$ is determined
by only ${\bf H}$ and ${\bf r}^m$ through some function $F$, assuming that $n$ is constant.
From (\ref{eq.relax}), it is obvious that $F$ can not be a polynomial of ${\bf H}$, therefore
RC is not a Krylov subspace method.

We now define the GRC method as
\eq
{\bf \Psi}^m_{\scriptscriptstyle \rm GRC} = G({\bf \phi}^{m-1},{\bf H},{\bf r}^m)
\qe
where $G$ is a linear operator.
In this study, we adopt that used in the CPRC method,
\eq
G({\bf \phi}^{m-1},{\bf H},{\bf r}^m) = ({\bf 1 - H}) \phi^{m-1} + {\bf r}^m
\qe
Then (\ref{eq.cprcpsi})
is considered to be a special case of the GRC method.
As a new vector generation by the matrix-vector operation for (2.5)
is repeated iteratively in GRC, it is obvious that it constructs
a Krylov subspace.
Furthermore, the initial residual is
\eq
{\bf r}^0 = {\bf b-HU}^0
\qe

then

\eq \phi^0 = \alpha_1^0 {\bf \Psi}^0
= \alpha_1^0 {\bf r}^0
\qe
\[
\equiv Q^0({\bf H}) {\bf r}^0
\]
where $Q^n({\bf H})$ denotes an $n$-th order polynomial of ${\bf H}$ and
the superscript of $\alpha$ indicates the iteration step $m$.
then the next residual
\eq
{\bf r}^1= {\bf r}^0-{\bf H\phi}^0
\qe
\[ = {\bf r}^0-{\bf H} Q^0({\bf H}) {\bf r}^0\]
\[
\equiv P^{1}({\bf H}) {\bf r}^0
\]
where $P^n({\bf H})$ denotes another $n$-th order polynomial of ${\bf H}$.
Likewise, the next $\phi$ becomes

\eq
\phi^1 = \alpha_1^1 {\bf \Psi}^1 + \sum_{k=2}^{L} \alpha_k^1 \phi^{2-k}
= \alpha_1^1 {\bf \Psi}^1 + \alpha_2^1 \phi^0
\qe
\[
= \alpha_1^1 \{ {( \bf 1 - H}) \phi^{0} + {\bf r}^1 \} + \alpha_2^1 \phi^0
\]
\[
= \alpha_1^1 \{ {( \bf 1 - H}) Q^0({\bf H}) {\bf r}_0 + {\bf r}^1 \} + \alpha_2^1 Q^0({\bf H}) {\bf r}_0
\]
\[
\equiv Q^1({\bf H}) {\bf r}_0
\]

and the next residual becomes

\eq
{\bf r}^2 = {\bf r}^1-{\bf H\phi}^1
\qe

\[
= P^{1}({\bf H}) {\bf r}^0 - {\bf H} Q^1({\bf H}) {\bf r}_0
\]

\[
\equiv P^2({\bf H}) {\bf r}^0
\]

Thus vectors $\phi^0, {\bf r}^1$ and $\phi^1, {\bf r}^2$
are expressed by polynomials of ${\bf H}$ times ${\bf r}^0$.
Furthermore, for $m=3,4$ and $L=3$, the vectors will be

\eq
\phi^2 = \alpha_1^2 {\bf \Psi}^2 + \sum_{k=2}^{L} \alpha_k^2 \phi^{3-k}
\qe
\[
= \alpha_1^2 \{ ({\bf 1 - H}) \phi^{1} + {\bf r}^2 \} + \sum_{k=2}^{L} \alpha_k^2 \phi^{3-k}
\]

\[
= \alpha_1^2 \{ ({\bf 1 - H}) Q^{1}({\bf H}) {\bf r}_0 + P^{2}({\bf H}) {\bf r}_0\}
+ \sum_{k=2}^{L} \alpha_k^2 Q^{3-k}({\bf H}) {\bf r}_0 
\]

\[
= \alpha_1^2 \{ ({\bf 1 - H}) Q^{1}({\bf H}) {\bf r}_0 + P^{2}({\bf H}) {\bf r}_0\}
+ \alpha_2^2 Q^{1}({\bf H}) {\bf r}_0
+ \alpha_3^2 Q^{0}({\bf H}) {\bf r}_0 
\]

\[
\equiv Q^2({\bf H}) {\bf r}^0
\]

and the residual becomes

\eq
{\bf r}^{3} = {\bf r}^2-{\bf H\phi}^2
\qe
\[
\equiv P^{3}({\bf H}) {\bf r}^0
\]

\eq
\phi^3 = \alpha_1^3 {\bf \Psi}^3 + \sum_{k=2}^{L} \alpha_k^3 \phi^{4-k}
\qe
\[
= \alpha_1^3 \{ ({\bf 1 - H}) \phi^{2} + {\bf r}^3 \} + \sum_{k=2}^{L} \alpha_k^3 \phi^{4-k}
\]

\[
= \alpha_1^3 \{ ({\bf 1 - H}) Q^{2}({\bf H}) {\bf r}_0 + P^{3}({\bf H}) {\bf r}_0\}
+ \sum_{k=2}^{L} \alpha_k^3 Q^{4-k}({\bf H}) {\bf r}_0 
\]

\[
= \alpha_1^3 \{ ({\bf 1 - H}) Q^{2}({\bf H}) {\bf r}_0 + P^{3}({\bf H}) {\bf r}_0\}
+ \alpha_2^3 Q^{2}({\bf H}) {\bf r}_0
+ \alpha_3^3 Q^{1}({\bf H}) {\bf r}_0 
\]

\[
\equiv Q^3({\bf H}) {\bf r}^0
\]

and the residual becomes

\eq
{\bf r}^{4} = {\bf r}^3-{\bf H\phi}^3
\qe
\[
\equiv P^{4}({\bf H}) {\bf r}^0
\]

Then for the consecutive residuals (and also $\phi^m$),
suppose
\eq{\bf r}^n= P^n({\bf H}) {\bf r}^0
\;\; {\rm and} \;\;
\phi^{n-1}= Q^{n-1}({\bf H}) {\bf r}^0
\;\; {\rm for} \;\;
n \leq m.
\qe

then

\eq
\phi^m = \alpha_1^m {\bf \Psi}^m + \sum_{k=2}^{L} \alpha_k^m \phi^{m-k+1}
\qe
\[
= \alpha_1^m \{ ({\bf 1 - H}) \phi^{m-1} + {\bf r}^m \} + \sum_{k=2}^{L} \alpha_k^m \phi^{m-k+1}
\]

\[
= \alpha_1^m \{ ({\bf 1 - H}) Q^{m-1}({\bf H}) {\bf r}_0 + P^{m}({\bf H}) {\bf r}_0\}
+ \sum_{k=2}^{L} \alpha_k^m Q^{m-k+1}({\bf H}) {\bf r}_0 
\]
considering $m-k+1$ in the summation is always less than $m$, we obtain
\[
\equiv Q^m({\bf H}) {\bf r}^0
\]

and the residual becomes

\eq
{\bf r}^{m+1} = {\bf r}^m-{\bf H\phi}^m
\qe
\[
\equiv P^{m+1}({\bf H}) {\bf r}^0
\]

Since $\phi^0, {\bf r}^1$ and $\phi^1, {\bf r}^2$
are expressed by polynomials of ${\bf H}$ times ${\bf r}^0$,
set $n=3$ and apply the above recurrence, namely, 
$\phi^2, {\bf r}^3$, $\phi^3, {\bf r}^4$ and so on are all
expressed by a polynomial of ${\bf H}$ times  ${\bf r}^0$,
As a result, the updated residual ${\bf r}^{m+1}$ for general $m$ is expressed
by $P^{m+1}({\bf H}){\bf r}^0$.

\eq
\phi^m = \alpha_1^m {\bf \Psi}^m + \sum_{k=2}^{L} \alpha_k^m \phi^{m-k+1}
\qe
\[
= \alpha_1^m \{ ({\bf 1 - H}) \phi^{m-1} + {\bf r}^m \} + \sum_{k=2}^{L} \alpha_k^m \phi^{m-k+1}
\]

If we assume $L>m+1$,

\[
= \alpha_1^m \{ ({\bf 1 - H}) \phi^{m-1} + {\bf r}^m \}
+ \sum_{n=0}^{m-1} \alpha_{m-n+1}^m \phi^n
\]
\[
= -\alpha_1^m {\bf H} \phi^{m-1} + \alpha_1^m {\bf r}^m
+ (\alpha_1^m+ \alpha_2^m) \phi^{m-1}
+ \sum_{n=0}^{m-2} \alpha_{m-n+1}^m \phi^n
\]

Therefore GRC is shown to be  a Krylov subspace method. 
Then, it is considered that its application is not limited
to coupled perturbed equations and that
it can be applied to general linear systems.
This will be an advantage to the RC method, whose application is limited
to linear problems for which a relaxation method is valid.

Since the subspace is limited in $L$ dimensions,
GRC is a truncated Krylov subspace method, while the GMRES(K) method
is a restarted method, where $K$ is the restart number.
Therefore GRC can be expressed by recurrence with only $L$ terms,
(typically $L$ is from 3 to 10), at the expense of the global orthogonality
not being guaranteed.

\subsection{Relation to the conjugate residual (CR)}

We have shown that $\phi^m$ is
expressed by polynomials of ${\bf H}$ times ${\bf r}^0$ as
\eq
\phi^m = Q^m({\bf H}) {\bf r}^0
\qe

which means

\eq
\phi^m \in {\cal K}_m
\qe

where

\eq
{\cal K}_i = {\cal K}_i({\bf H},{\bf r}_0) = {\rm span}\{{\bf r}_0,{\bf H}{\bf r}_0,{\bf H}^2{\bf r}_0,\cdots,{\bf H}^{i-1}{\bf r}_0 \}
\qe

is a Krylov subspace of degree $i$.
Moreover,

\eq
\Psi^m \in {\cal K}_m
\qe

The new direction vector $\phi^m$ is generated so that
it becomes ${\bf H}^T{\bf H}$-orthogonal to
previous ones,

\eq
({\bf H} \phi^m , {\bf H} \phi^j) = 0, j<m
\qe

therefore

\eq
{\bf H} \phi^m \bot {\bf H} {\cal K}_{m-1}
\qe

This is also equivalent to Gram-Schmidt  ${\bf H}^T{\bf H}$-orthogonalization
\eq
 \phi^{m}  = \Psi^m
- \sum_{n=0}^{m-1} ({\bf H} \phi^m , {\bf H} \phi^n) \phi^n
\qe

In CR, the direction vector corresponding to $\phi^m$ in GRC is
${\bf p}^m$, except that $\phi^m$ is the update for the current
solution but the update in CR is $a_m{\bf p}^m$ where
\eq
a_m = ({\bf H r}_m, {\bf r}_m) / ({\bf H p}_m,{\bf H p}_m)
\qe
which minimizes the updated residual norm.
Just like GRC, direction vectors are generated by
the polynomial of $\bf H$ times ${\bf r}_0$, therefore

\eq
{\bf p}^m \in {\cal K}_m
\qe

\eq
{\bf H} {\bf p}^m \bot {\bf H} {\cal K}_{m-1}
\qe

Moreover, both  $\phi^m$ and $a_m{\bf p}^m$ are determined
so that ${\bf r}^{m+1} = {\bf r}^m - {\bf H} \phi^m$
and ${\bf r}^{m+1} = {\bf r}^m - {\bf H} a_m \phi^m$
correspondingly are minimized.
As a result, if $\bf H$ is symmetric, which is a necessary condition for CR,
not only the Krylov subspaces but also the direction vectors
(with the scalar multiplied for CR) are identical
\eq
\phi^m  = a_m{\bf p}^m
\qe
for each iteration.
In the above discussion, the $L$ is assumed to be large enough
to guarantee ${\bf H}^T{\bf H}$-orthogonalization for all $\phi_m$.
In fact, for the symmetric $\bf H$ it appears that $L=3$ is sufficient
by applying the  ${\bf H}^T{\bf H}$-orthogonality relationship among the
direction vectors
and the residual vectors, and also the conjugate relationship among the residual vectors,
which is the theoretical result from CR.

Furthemore, if we let
\eq \label{eq.psiresid}
{\bf \Psi}^m = {\bf r}^m
\qe
then it is completely the same procedure in the CR iteration and
$L=2$ is sufficient.
It may be possible that this simpler expression can be used instead of the current $\Psi$ in (\ref{eq.cprcpsi}).
Experimental results are shown in
Figures \textcolor{green!50!black}{1} - \textcolor{green!50!black}{3}
to show the equivalence to CR and also
the difference between using eqs. (\ref{eq.cprcpsi}) and (\ref{eq.psiresid}).

\begin{figure}[!h]

\bgc
\caption{
Results of applying CR and GRC(L=10) to a symmetric matrix,
which originates from 'olm100' of the University of Florida database.
It is obtained by averaging the matrix and its transpose.
}
\edc
\includegraphics[scale=0.38]{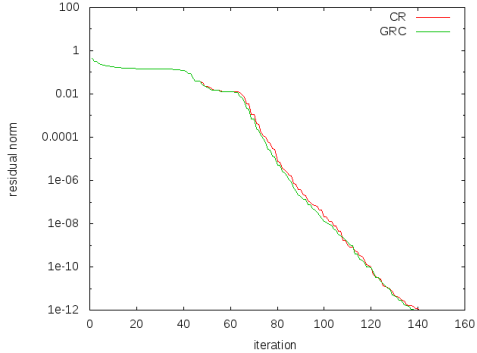}
\label{fig2.1}

\end{figure}

\begin{figure}[!h]
\bgc
\caption{
Results of applying GRC to the symmetric matrices of
Figure \textcolor{green!50!black}{1} and the matrix (1) in the Sec. 3.2
(which is also turned out to be symmetric).
'GRC2' denotes GRC(L=2) and so on.
GRC with $L \geq 3$ generates almost identical residual norms
at the beginning
in both graphs, as predicted by theory.
The difference at later iterations in the left
seems due to the accumulation of floating point errors.
}
\edc
\includegraphics[scale=0.38]{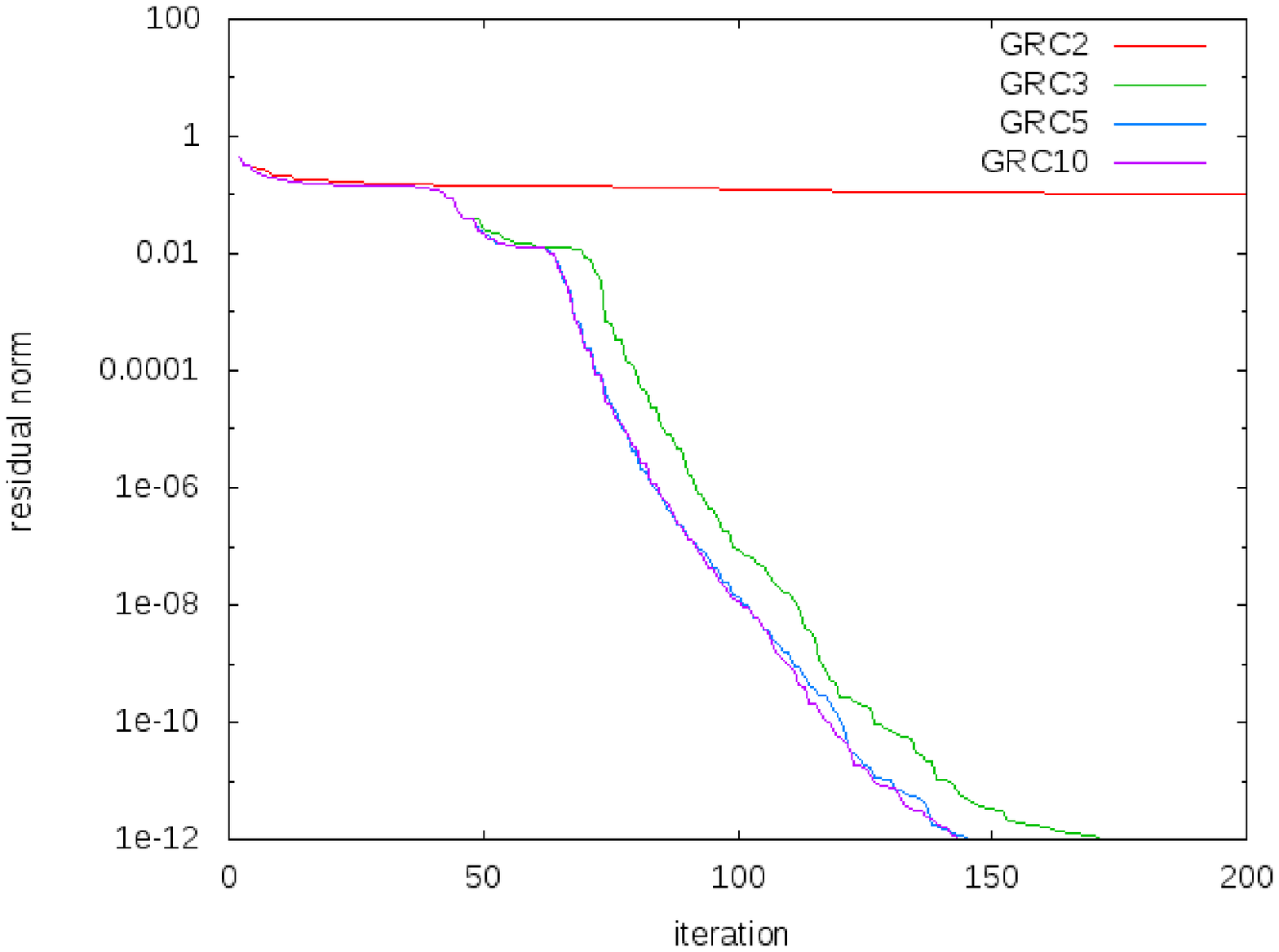}
\includegraphics[scale=0.38]{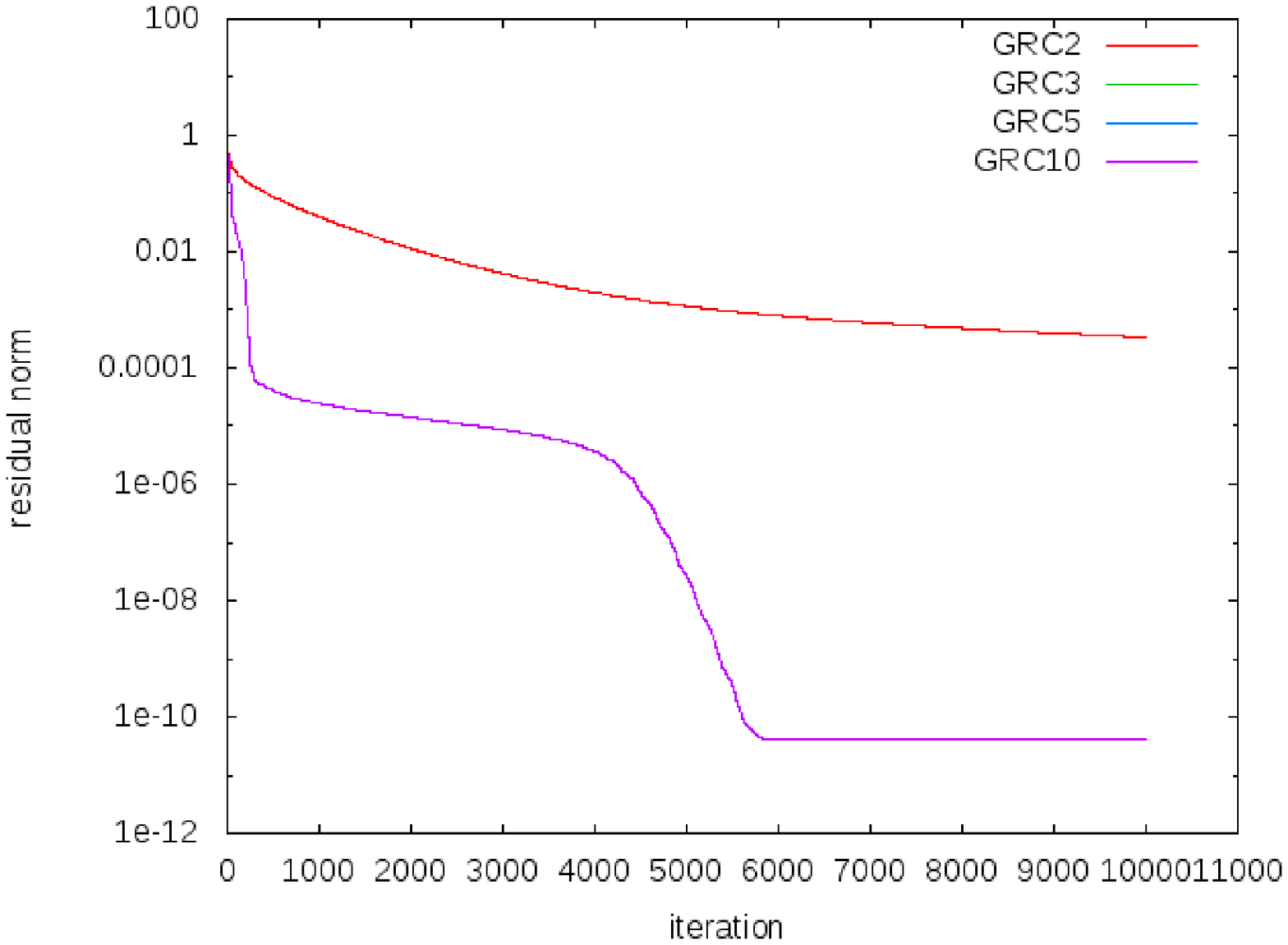}
\label{fig2.2}

\end{figure}

\begin{figure}[!h]
  
\bgc

\vspace{3mm}

\caption{
Results of applying GRC using eq. (\ref{eq.psiresid})
instead of eq. (\ref{eq.cprcpsi})
to the symmetric matrices of
Figure \textcolor{green!50!black}{1} and the matrix (1) in the Sec. 3.2.
In this case,  $L \geq 2$ is sufficient in theory.
In the right, all the contours overlap completely.
}
\edc

\includegraphics[scale=0.38]{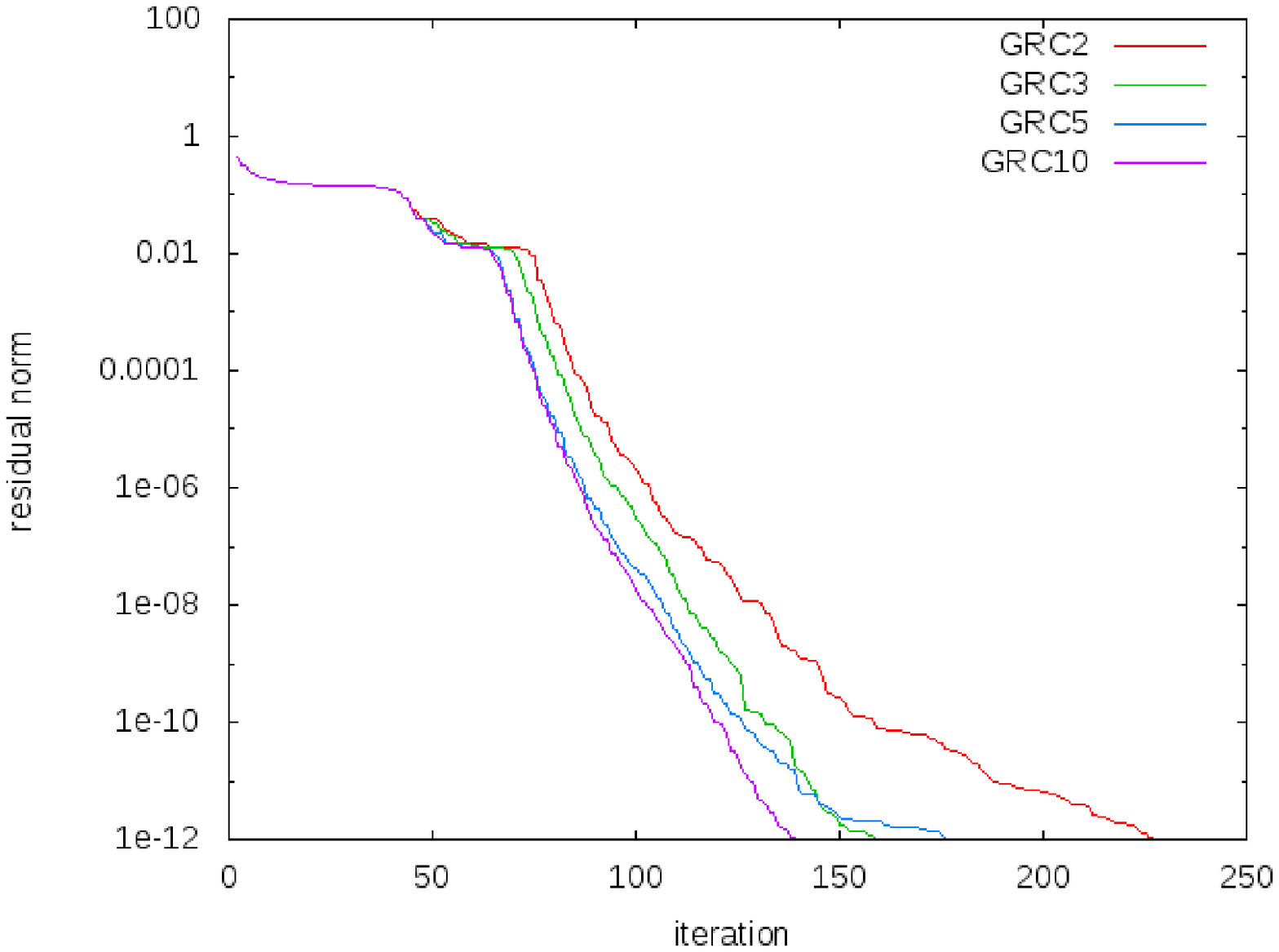}
\includegraphics[scale=0.38]{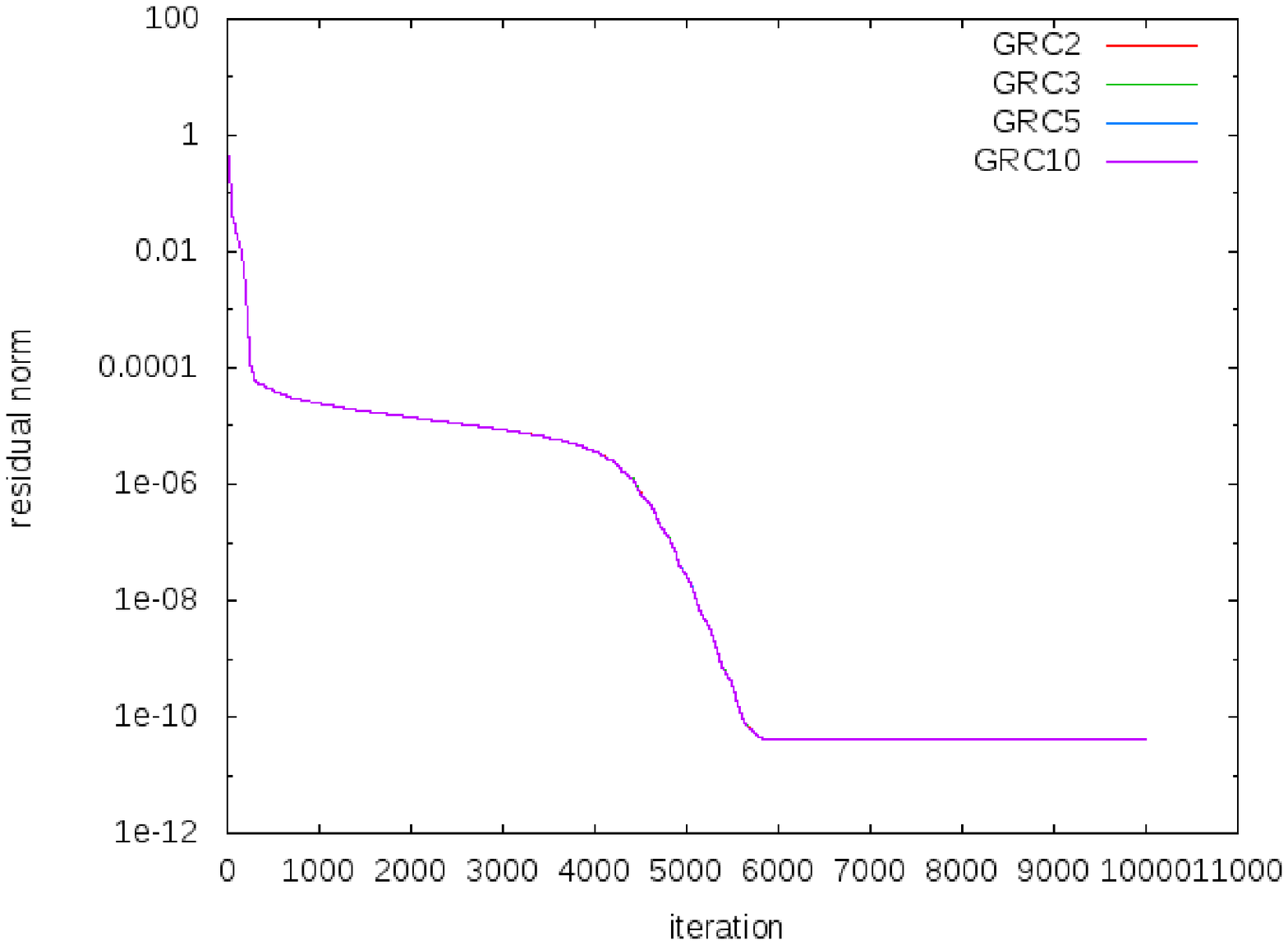}
\label{fig2.3}

\end{figure}

%\newpage

\section{Numerical experiments}

We evaluate the performance of GRC as compared to the original RC,
BiCGSTAB and GMRES,
with popular coefficient matrices used in other references.
Also, we investigate relationship between the size and convergence
with matrices whose sizes are determined by a parameter.

\subsection{Experimental condition}
For our experiments, we use
BiCGSTAB and GMRES that are implemented as subroutines
in the LIS library\cite{ref_lis}.
For this purpose, we have implemented the GRC and RC methods as LIS subroutines
in the same manner as BiCGSTAB and GMRES, so that
they use the same LIS routines such as matrix and vector operations.
As for the dimension of the subspace for both of GRC and RC, we set $L=5$ 
%In order to avoid influence of numerical errors when the vector size is large,
%inner product of the pair of basis vectors is calculated
%also when the orthogonality condition holds therefore theoretical inner product becomes zero.
%In that case, the inner product becomes a very small number,
%which will be used as a component of the coefficient matrix for the residual minimization equation.
As the inner solver of the RC method, SOR is used with the relaxation coefficient and
the number of iterations is set equal to 1.9 and 50, respectively.
The restart number for GMRES is set to be 40, which is the default value with LIS.
Calculation is done on a single core of Intel Core2 (3GHz) processor with 8 GB memory.

\subsection{Result with the test matrices}
The following three matrices are used for evaluation.
The first one is the coefficient matrix that we have been using for evaluating the RC method.
The other two are popular coefficient matrices among references
(for example \cite{ref_gmres}\cite{ref_bicgstab}).

\noindent
(1) Coefficient matrix generated by discretizing a Poisson equation on non-uniform grid
with a Neumann boundary condition,
which makes it difficult for relaxation methods to converge.

\noindent
(2) Coefficient matrix generated by discretizing the partial differential
equation \cite{ref_pdematrix}
\[
u_{xx} + u_{yy} + u_{zz} + 1000u_x = F
\]

\noindent
(3) Coefficient matrix named 'raefsky2' from the database of University of Florida sparse matrix collection
\cite{ref_florida}.
The right hand side vector for the linear system is set so that the solution vector
will be the vector with all the elements being unity.

Figures \textcolor{green!50!black}{4} - \textcolor{green!50!black}{6} show the results with these matrices.
The residual norm versus iteration step is shown in the left
and the residual norm versus time in the right.

With matrix (1),  RC resulted in the fastest convergence in time.
This is presumably because the relaxation method as its inner solver
converges efficiently with the discretized Poisson equation.
On the other hand, GRC shows slow convergence at the beginning
and accelerates gradually in later iteration steps. BiCGSTAB shows
specific fluctuation due to the lack of monotonically decreasing in residual norm.
For GMRES, residual norm stopped decreasing before convergence.
It should be noted that for RC, elapsed time in
a single iteration step is much longer than other methods, because
there are many inner iterations (50 in this case) of SOR in its inner solver.
For this reason, RC takes much less iteration steps compared to the elapsed time.

With matrix (2), the residual norm with the RC method does not decrease at all.
This is because SOR in the inner solver
diverged. Since the residual norm with the RC method in principle does not increase,
it remained constant.
%In addition, SOR converged by modifying its relaxation coefficient from the default value of 1.9
%to a smaller value of 1.0.
%GRC and GMRES show similar contours of convergence.

With matrix (3), the inner solver in the RC method also diverges
and the residual 
norm remained constant until some point, then it diverged.
This is due to the excessive degree of divergence by the inner solver,
with the residual norm in the order of $10^{100}$ and residual minimization
does not work precisely any more. BiCGSTAB shows significantly fast
convergence. GMRES shows slower convergence, and GRC shows even slower convergence.

\begin{figure}[!t]

\bgc
\caption{Results from the coefficient matrix (1).}
\edc
\includegraphics[scale=0.28]{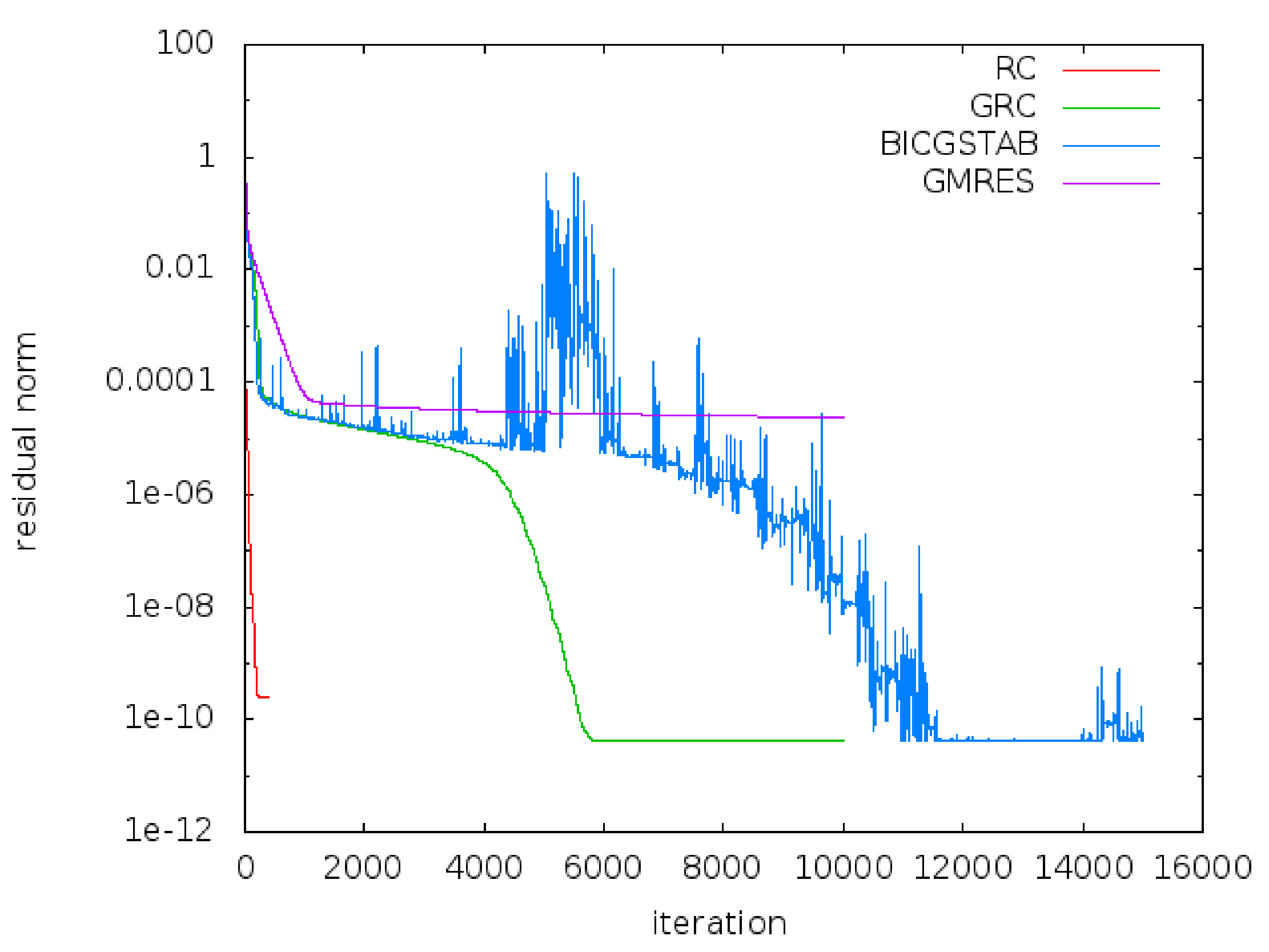}
\includegraphics[scale=0.28]{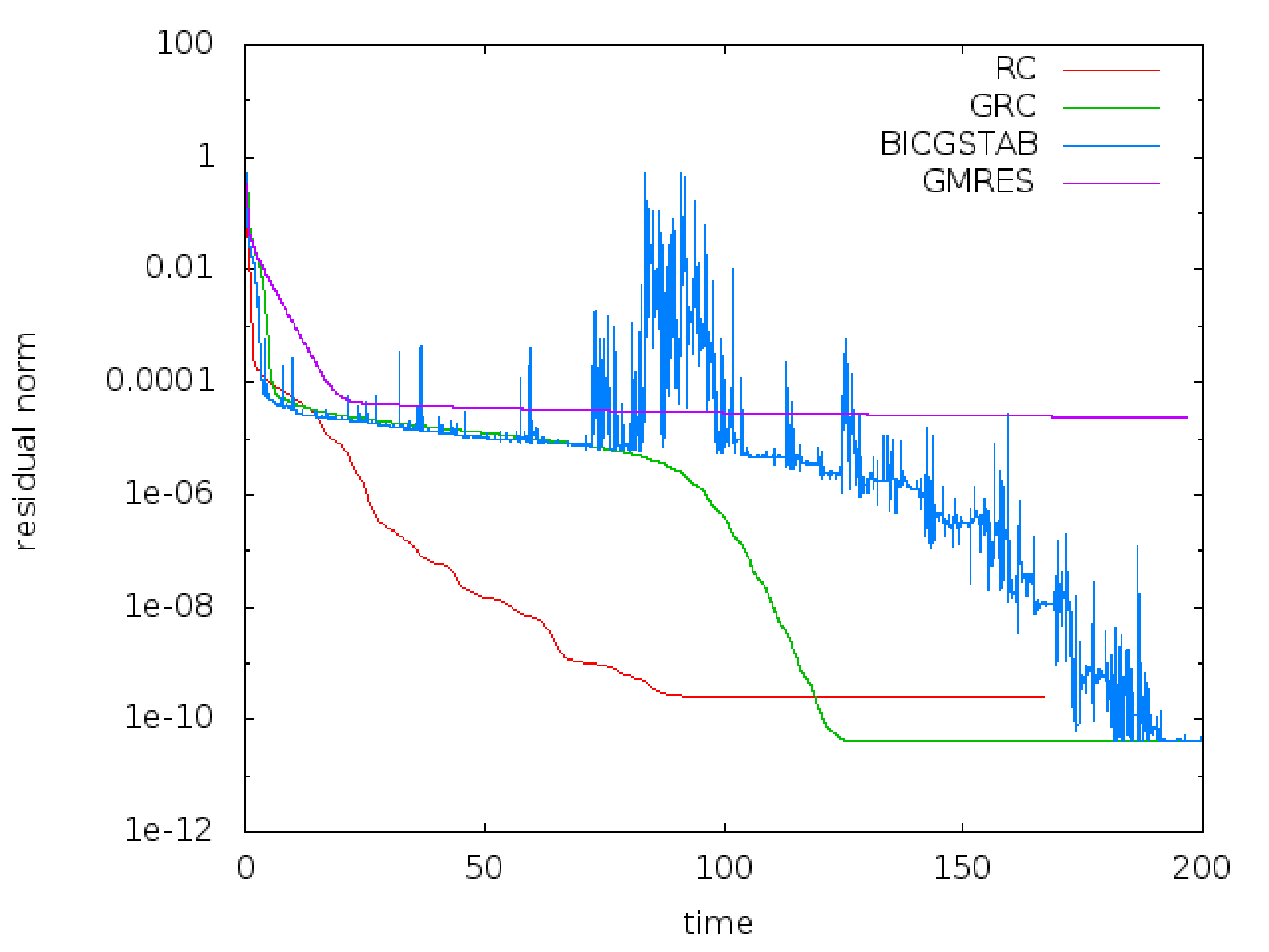}
\label{fig3.1}

\end{figure}
\begin{figure}[!t]
\bgc

\vspace{3mm}

\caption{Results from the coefficient matrix (2).}
\edc
\includegraphics[scale=0.28]{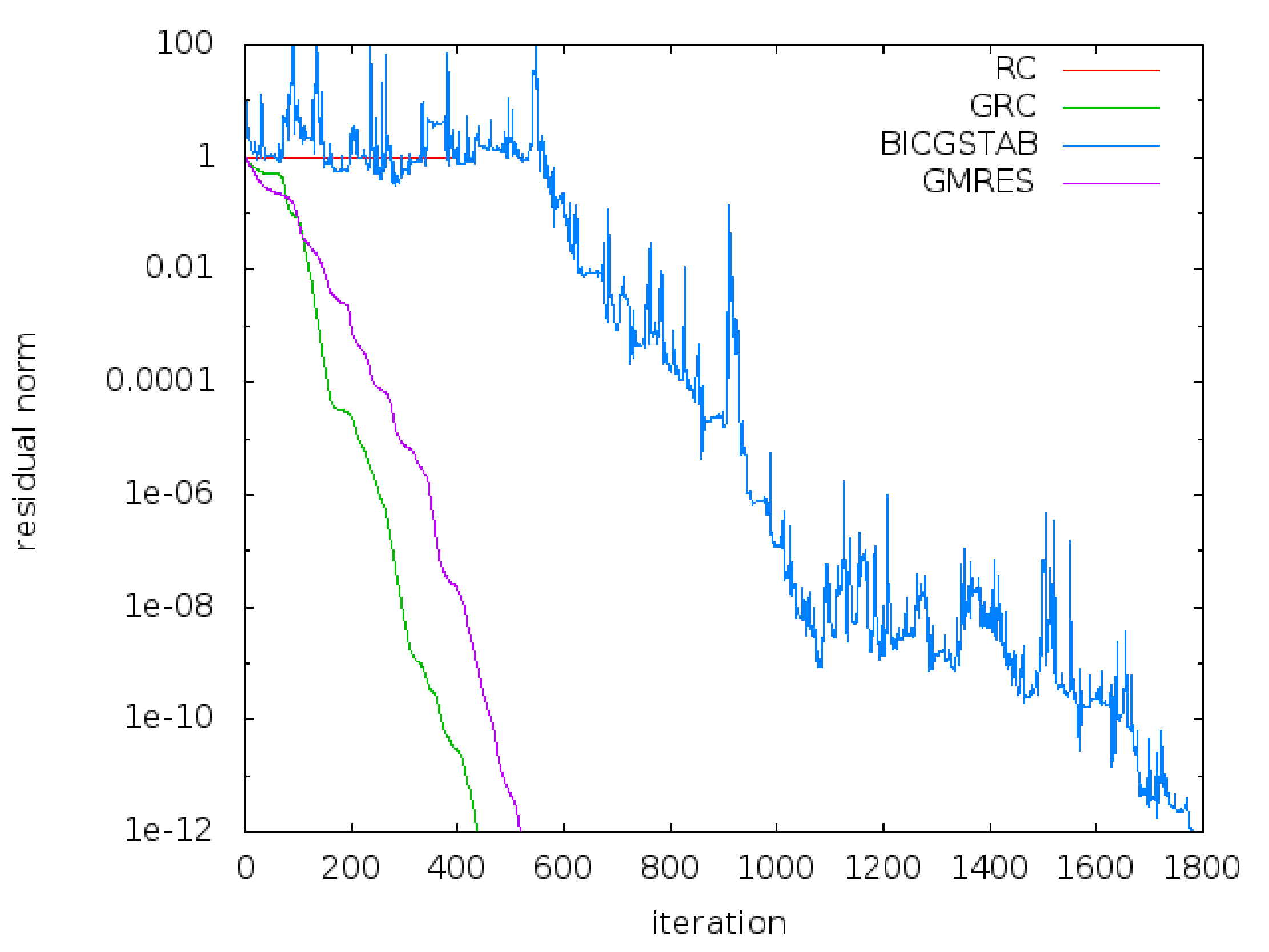}
\includegraphics[scale=0.28]{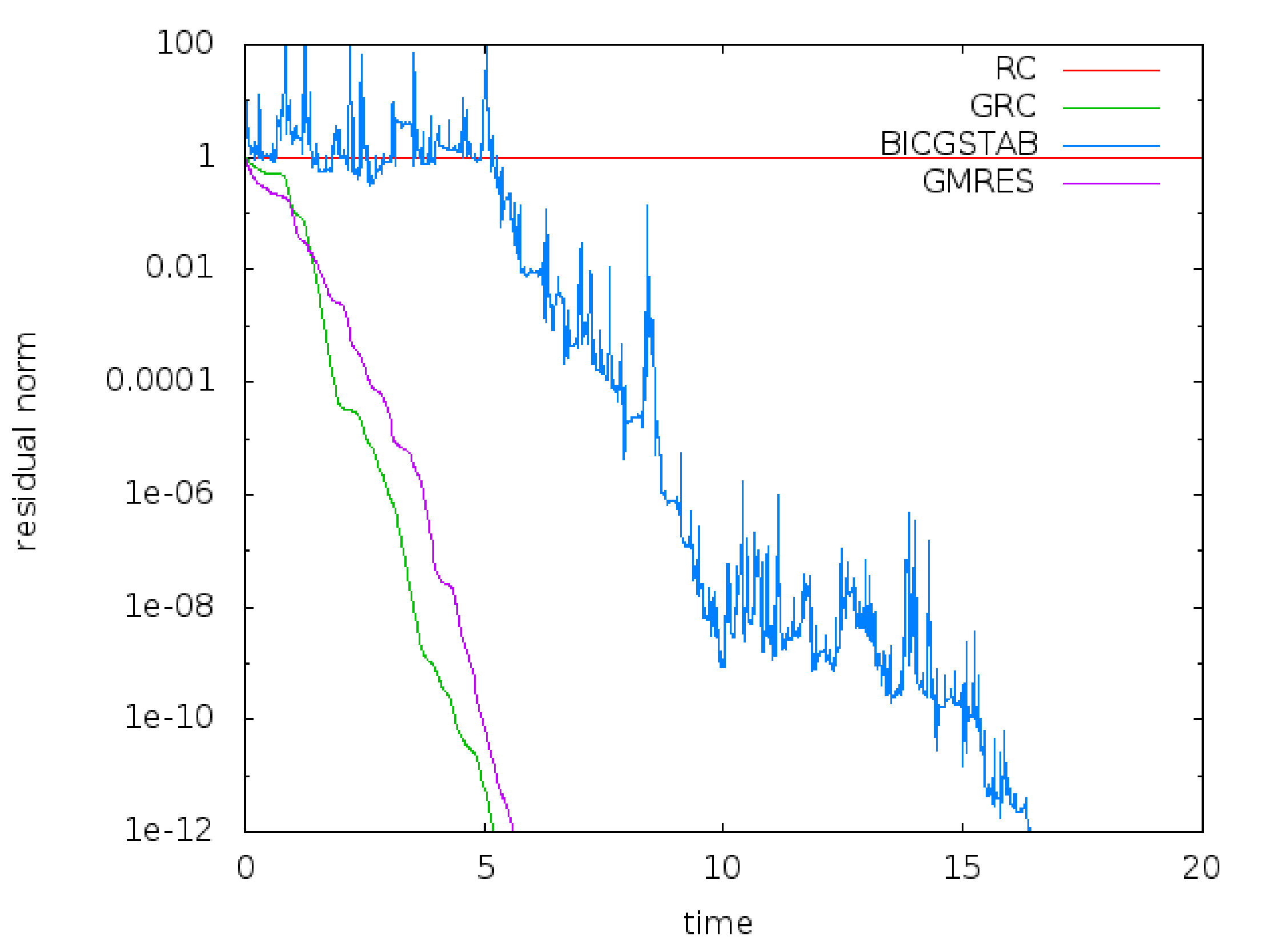}
\label{fig3.2}

\end{figure}
\begin{figure}[!t]
\bgc

\vspace{3mm}

\caption{Results from the coefficient matrix (3).}
\edc
\includegraphics[scale=0.28]{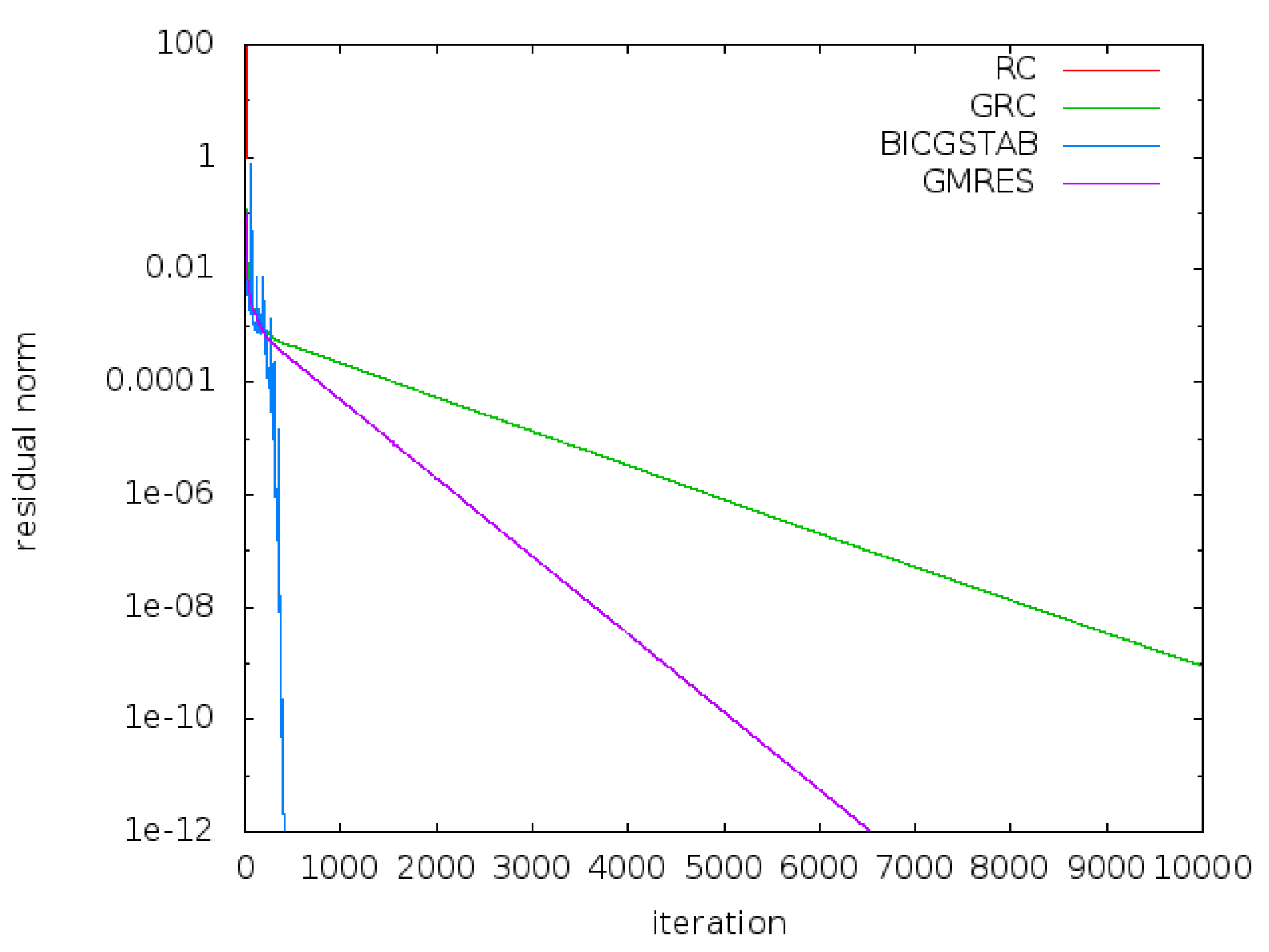}
\includegraphics[scale=0.28]{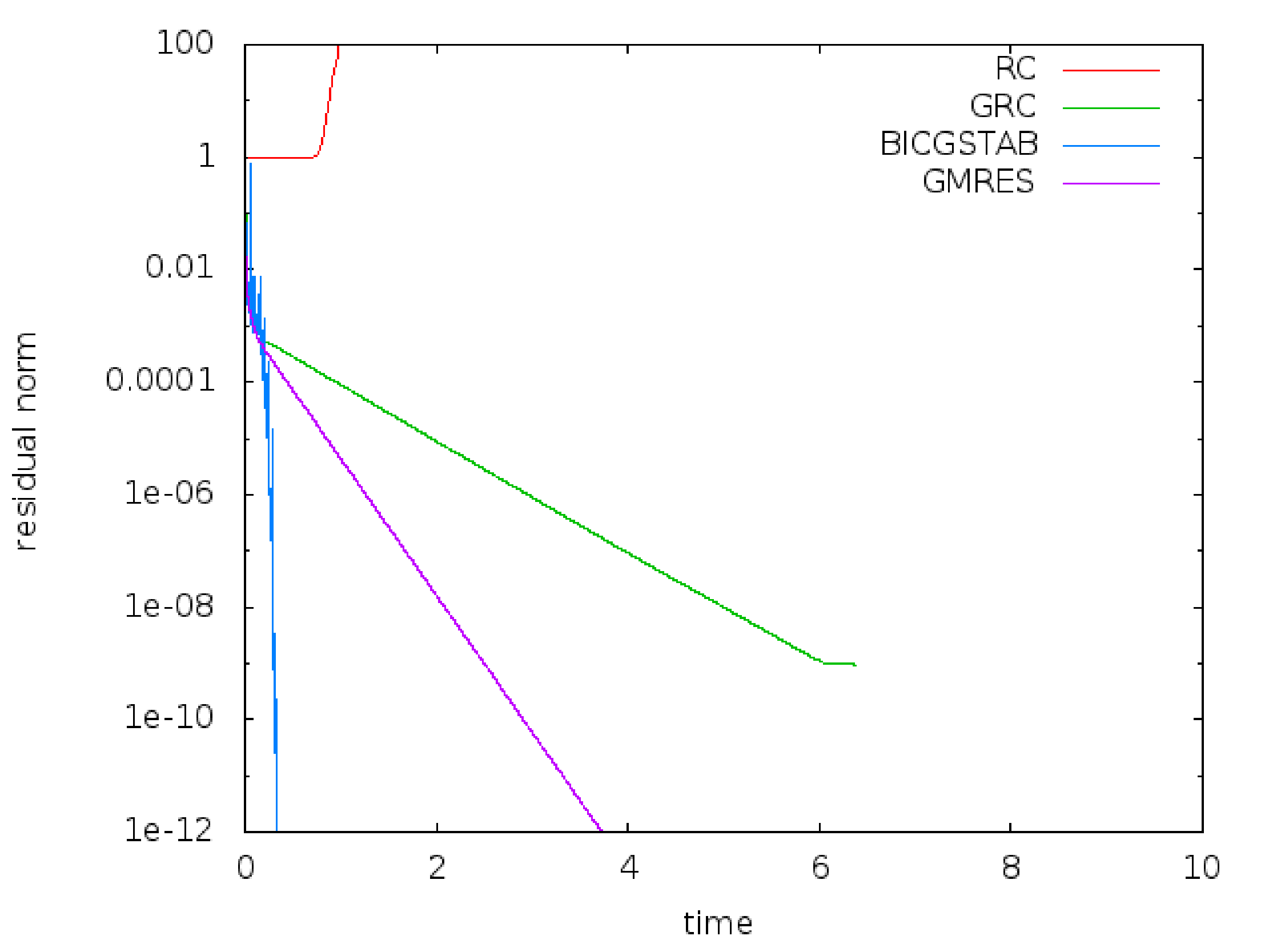}
\label{fig3.3}
%%%\includegrfaphics[width=0.8\textwidth,natwidth=610,natheight=642]{iter_standard.png}
\\
\\
\\

\end{figure}

%\setlength\fboxsep{0pt}
%\setlength\fboxrule{0.5pt}
%\fbox{\includegraphics{iter_standard.png}}

\begin{figure}[!h]
\bgc
\caption{Results from the coefficient matrix (2a).}
\edc
\includegraphics[scale=0.28]{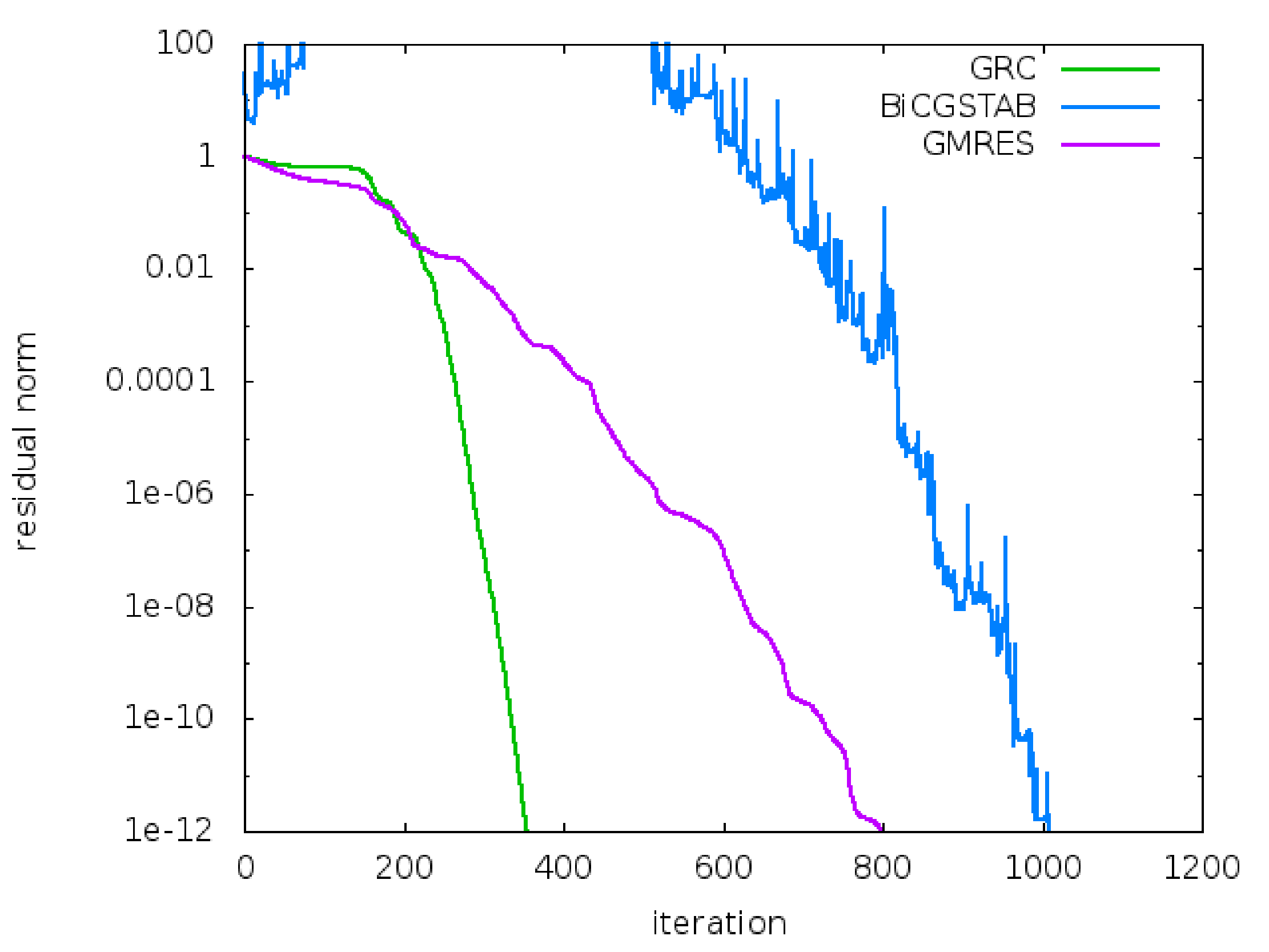}
\includegraphics[scale=0.28]{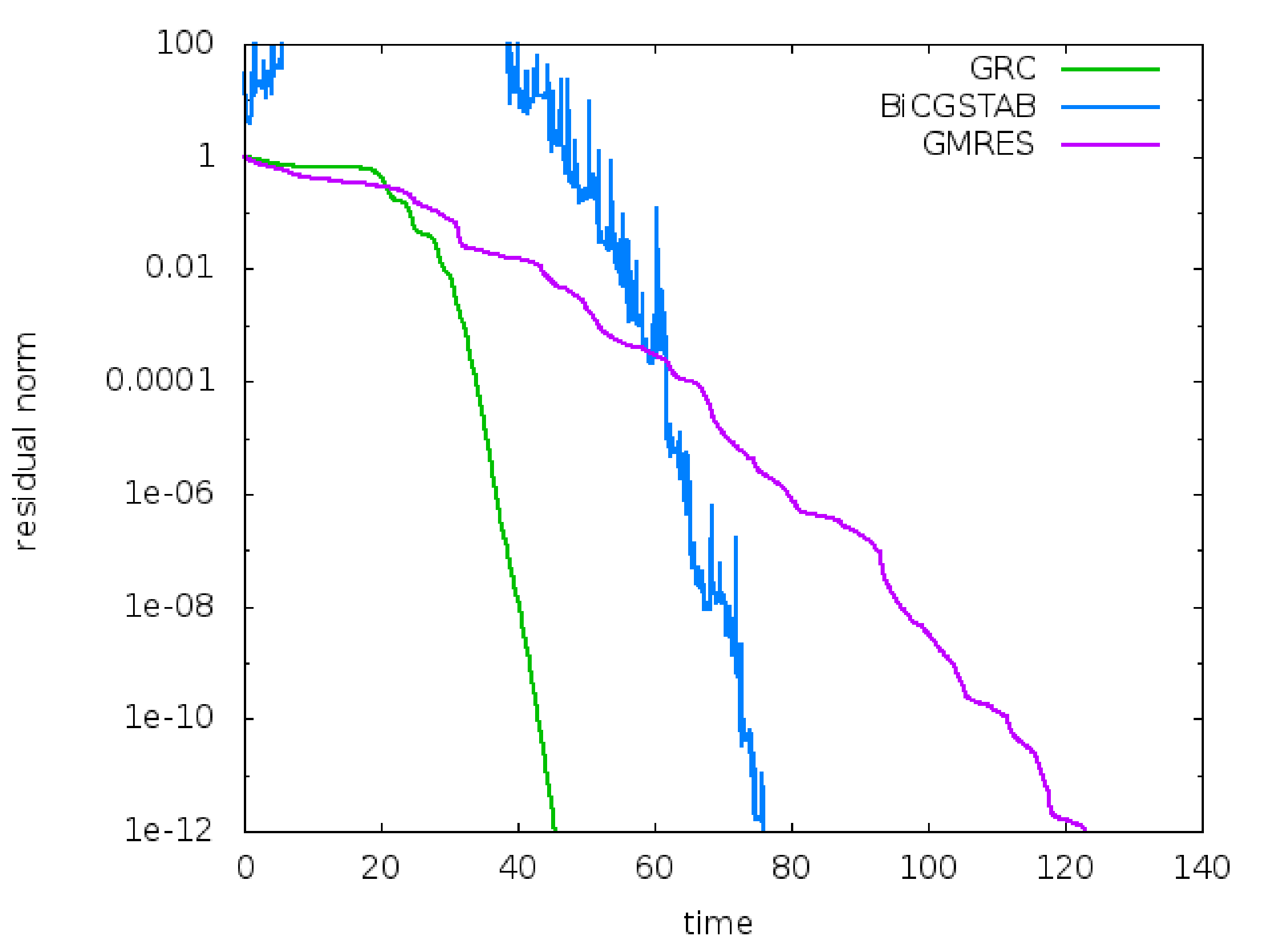}
\label{fig3.4}

\end{figure}
\begin{figure}[!h]

\bgc

\vspace{3mm}

\caption{Results from the coefficient matrix (2b).}
\edc
\includegraphics[scale=0.28]{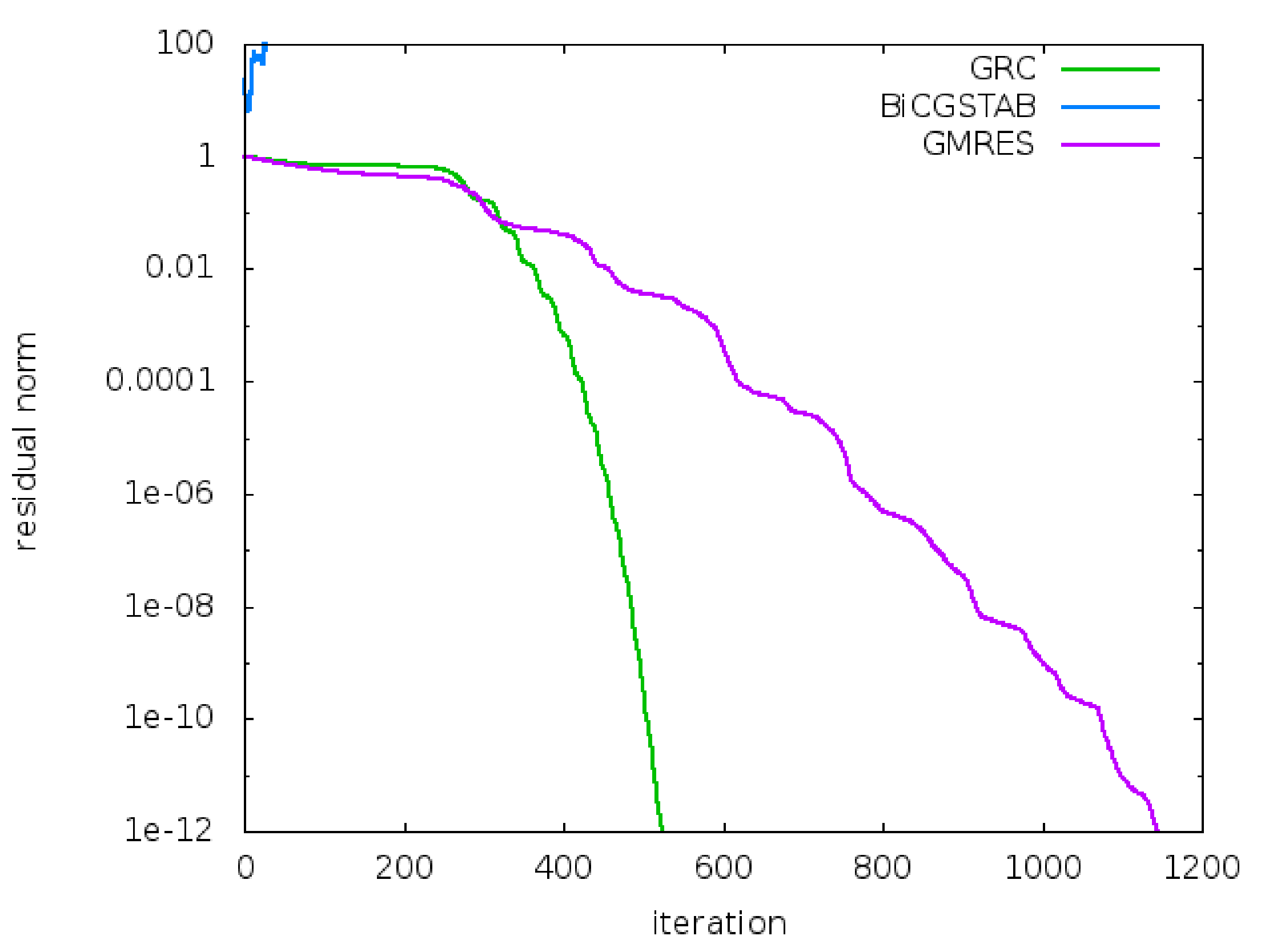}
\includegraphics[scale=0.28]{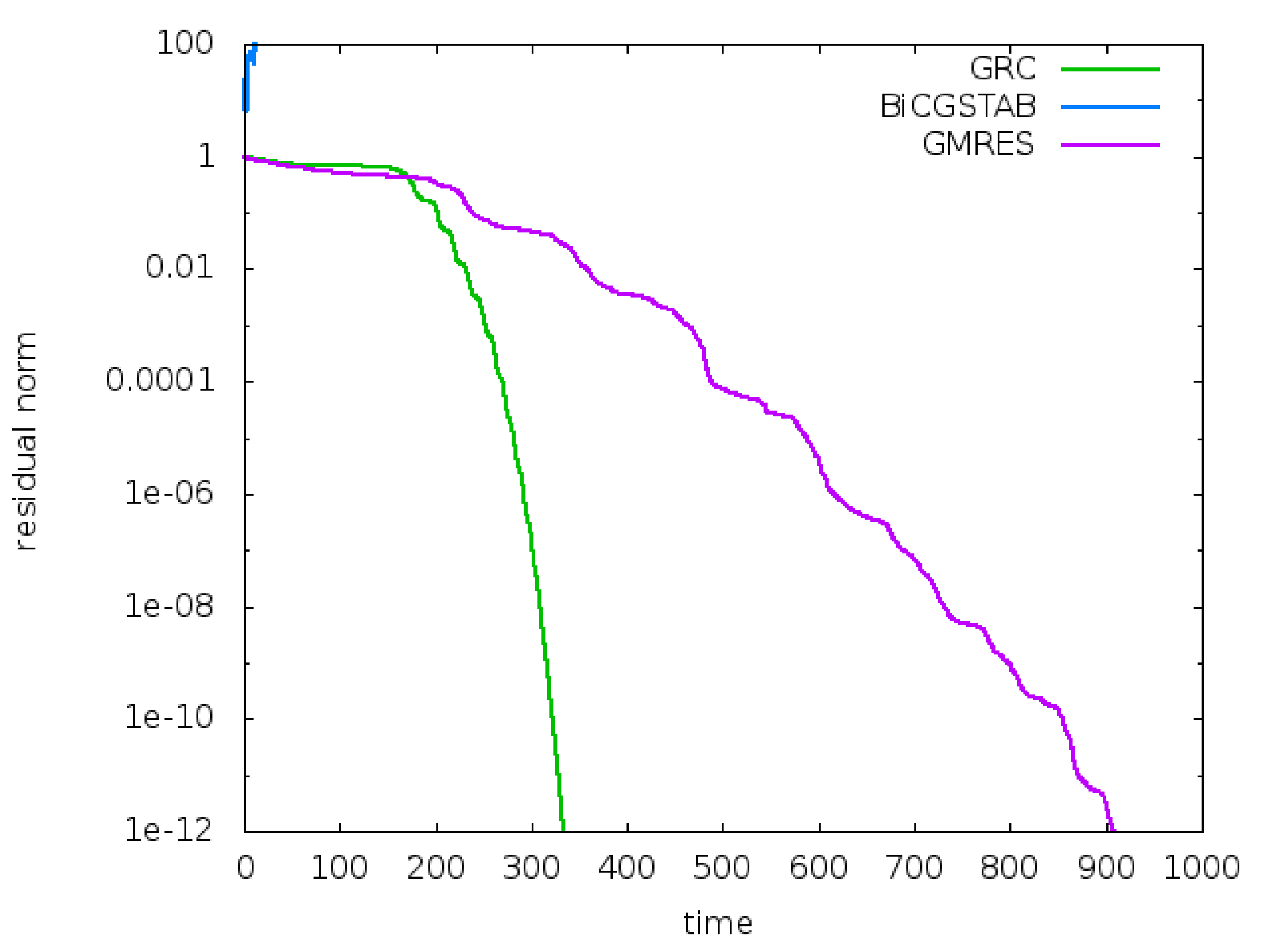}
\label{fig3.5}
\end{figure}

\begin{figure}[!h]
\bgc
\caption{Results from the coefficient matrix (2c).}
\edc
\includegraphics[scale=0.28]{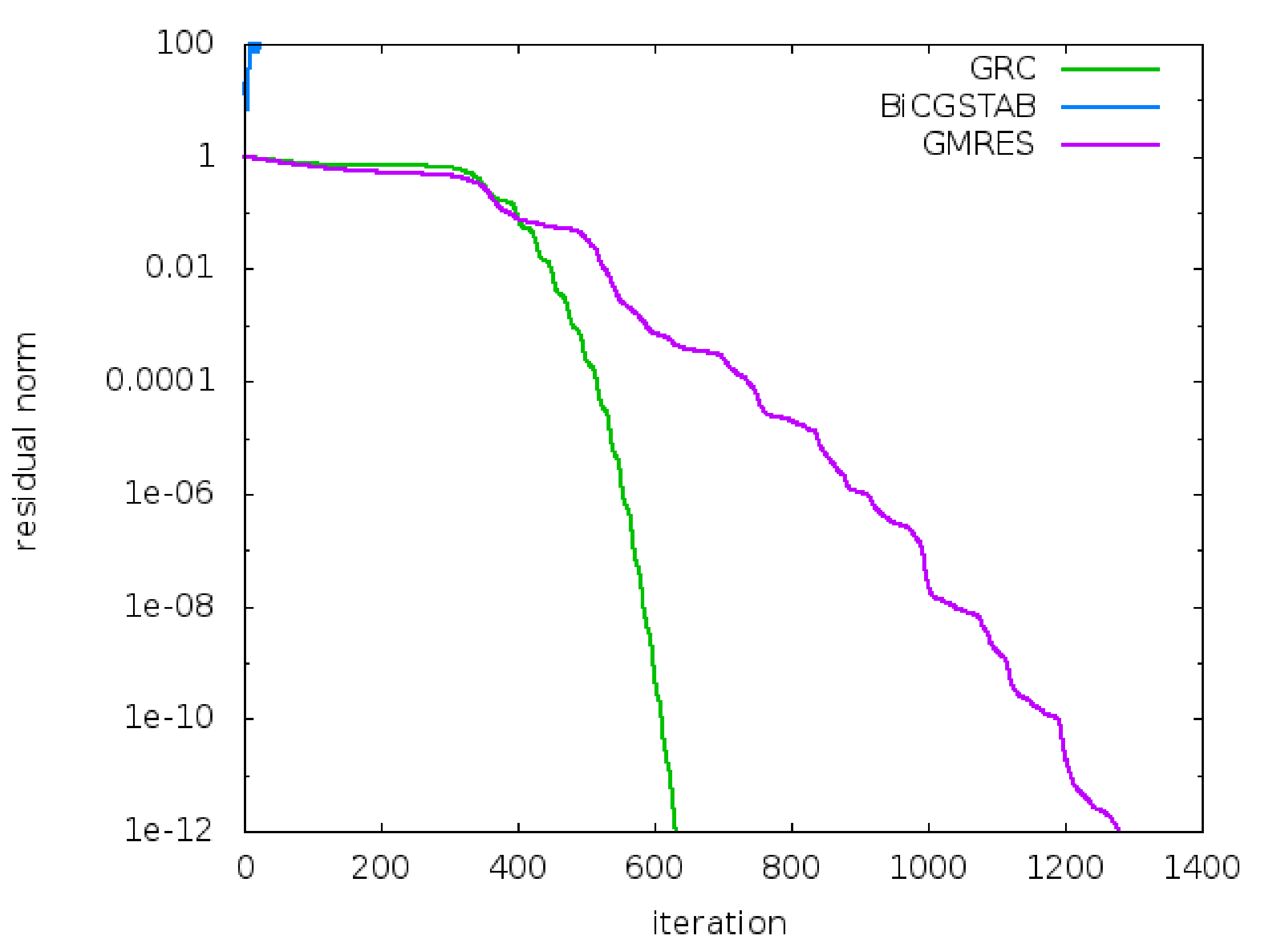}
\includegraphics[scale=0.28]{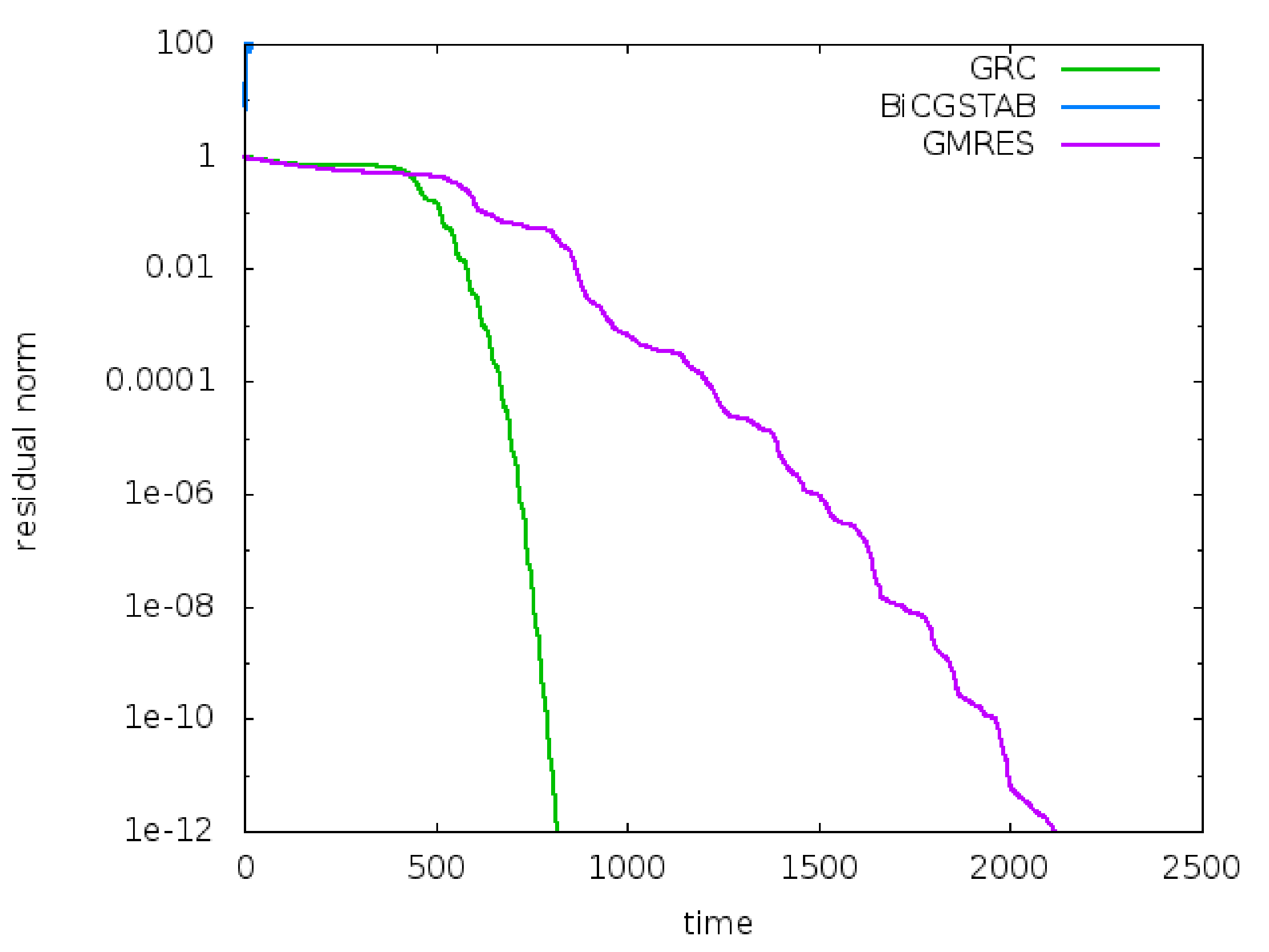}
\label{fig3.6}

\end{figure}
\begin{figure}[!h]
  
\bgc

\vspace{3mm}

\caption{Results from the coefficient matrix (2d).}
\edc
\includegraphics[scale=0.28]{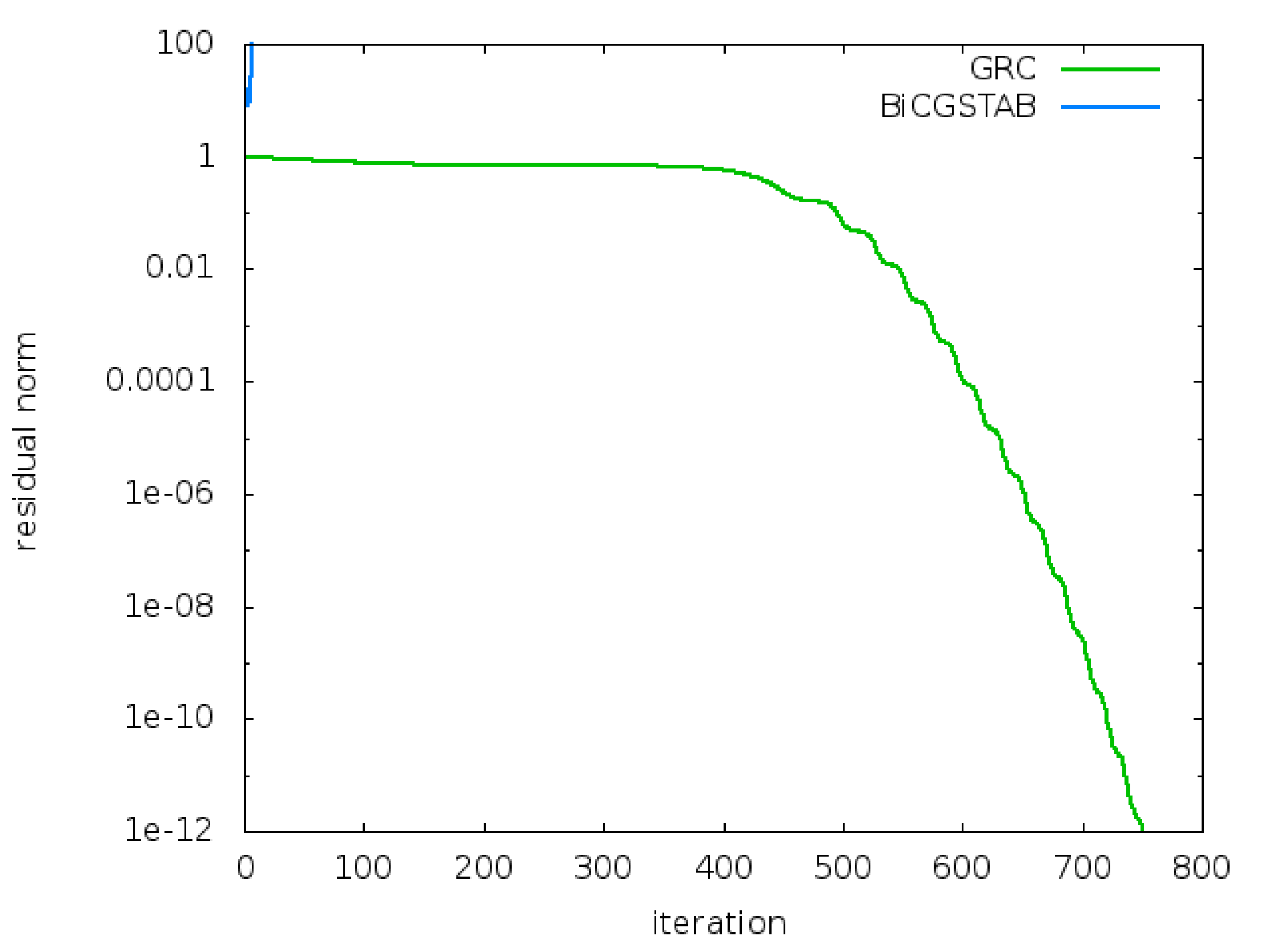}
\includegraphics[scale=0.28]{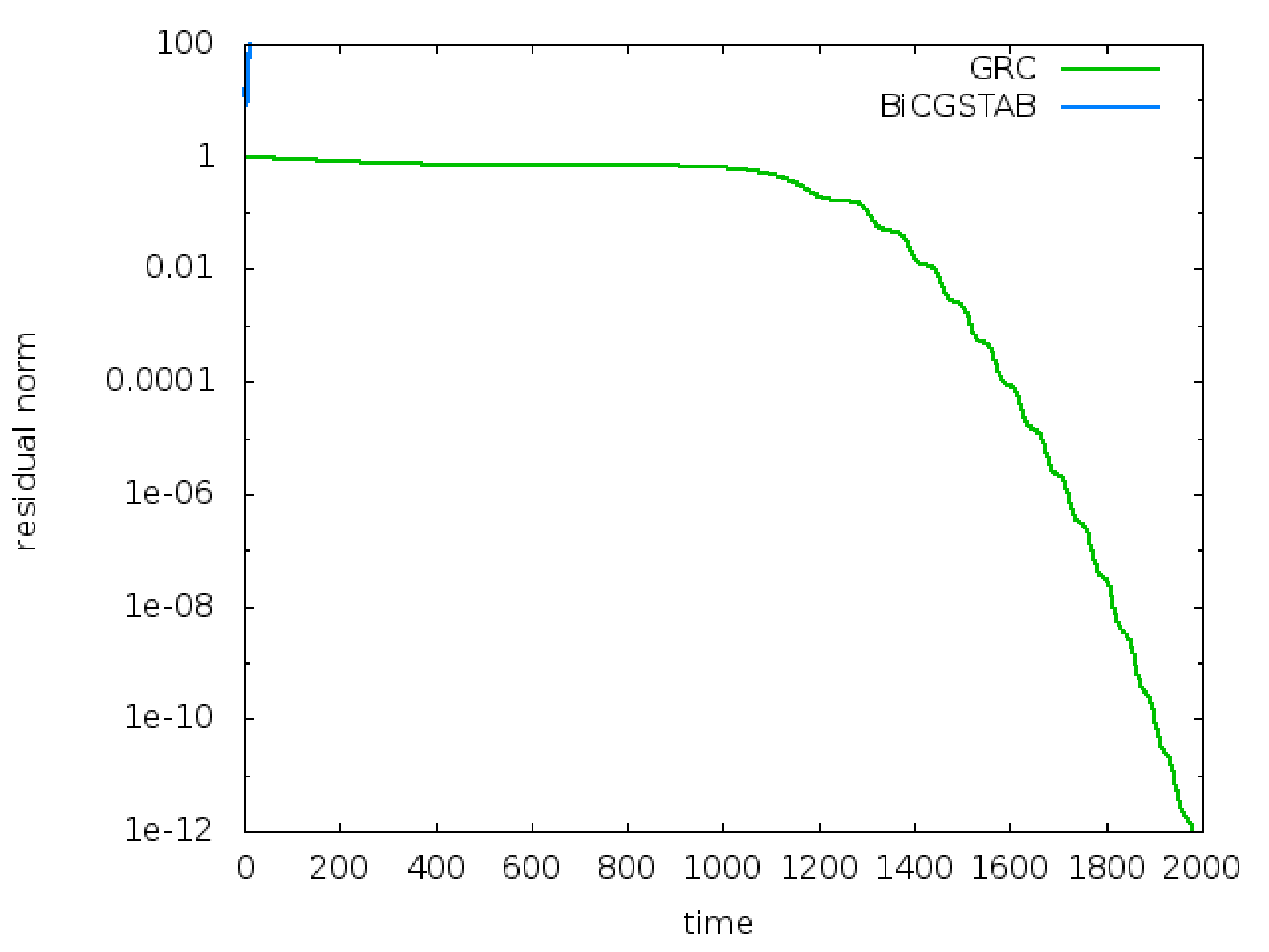}
\label{fig3.7}
\end{figure}

%\setlength\fboxsep{0pt}
%\setlength\fboxrule{0.5pt}
%\fbox{\includegraphics{iter_standard.png}}

\subsection{Dependence on size of the matrix from a partial differential
equation}
Size of the coefficient matrix (2) in the previous section can be
controlled by changing the fineness in discretization.
Thus we investigate convergence with the different size of the matrix.
Figures \textcolor{green!50!black}{7} - \textcolor{green!50!black}{10}
show the results for matrices shown in the table below.

\begin{table}[htbp]
\caption{Size of matrices.}
\begin{center}\footnotesize
\begin{tabular}{|c|c|c|c|c|}\hline
matrix    & 2a      & 2b      & 2c       & 2d        \\ \hline
size      & 1000000 & 4913000 & 9938375  & 19683000  \\ \hline
\#nonzeros & 6940000 & 34217600 & 69291275 & 137343600 \\ \hline
\end{tabular}
\end{center}
\end{table}

With the matrix (2a), the residual norm of BiCGSTAB becomes as large as about $10^4$,
then it converges. However, it diverges with larger matrices (2b), (2c) and (2d).
For GRC , both the number of steps and elapsed time until convergence are
less than half of those for GMRES.
In addition, with matrix (2c), GMRES terminated due to an out-of-memory error,
with 8GB memory of the experimental condition and as a result, GRC is the only
method that converged.

\subsection{Some other large matrices}
Figures \textcolor{green!50!black}{11} - \textcolor{green!50!black}{13}
show the results with large matrices from the aforementioned
University of Florida sparse matrix collection \\ database.
Unlike the previous section's examples, BiCGSTAB achieves best convergence with all these matrices.
Thus BiCGSTAB shows fast convergence if its residual norm does not diverge.
On the other hand, GRC and GMRES are robust in a sense that the residual norm
does not increase in principle, however, there are some cases where convergence stops on the way.
Also, as is the case with large size matrices in the previous section,
GRC tends to converge faster than the GMRES method.
One reason for this may be due to be the fact that the number of
basis vectors $L$ of GRC
is smaller than that of GMRES $K$ and therefore less sensitive to numerical errors
in residual minimization, which takes many inner product operations among basis vectors.

\section{Discussion}
The following table shows necessary memory (the unit is number of vectors), number of
matrix-vector multiplications (MATVEC), number of inner product operations (DOT)
in a single iteration step. If the restart number of GMRES, which is set $K=40$, is set smaller
such as $K=10$,
necessary memory naturally becomes smaller accordingly, however, characteristics on convergence
reportedly degrades significantly \cite{ref_gmres}. We also confirmed this fact in our preliminary
experiments.

\begin{table}[htbp]
\caption{Necessary memory, number of matrix-vector multiplications, and
number of inner product operations.
Constants $L$ and $K$ are the numbers of basis vectors for residual minimization
in GRC and GMRES, respectively.}
\begin{center}\footnotesize
\begin{tabular}{|c|c|c|c|}\hline
  & memory & MATVEC & DOT  \\ \hline
GRC      & $2L$ & 1 & $2L$ \\ \hline
BiCGSTAB & 5    & 2 & 2    \\ \hline
GMRES    & $K$  & 1 & $K$  \\ \hline
\end{tabular}
\end{center}
\end{table}

Among the four methods, RC shows the fastest convergence with
matrices for which the relaxation method of its inner solver converges,
for example, the test matrix (1).
It did not converge with other test matrices, although
it may converge with a much smaller relaxation coefficient.
The advantage of the RC method is that
it can further accelerate convergence if the inner solver converges.

Among the Krylov subspace methods (GRC, BiCGSTAB and GMRES), \\
BiCGSTAB needs to keep the least number of vectors.
GRC and GMRES, which guarantee that the residual norm does not increase in principle,
tend to show robust convergence, in contrast to BiCGSTAB which sometimes diverges.
Their contours of convergence sometimes look similar, probably due to the common principle
of residual minimization.
However, there are a few cases where GRC and GMRES fail to converge
while BiCGSTAB achieves fast convergence.
As for memory usage by GRC and GMRES,
GMRES needs to keep a relatively large number of vectors (restart number)
for effective convergence. On the other hand,
GRC $(L=5)$ needs to keep less number of vectors.

\section{Conclusion}

\begin{figure}[!h]
\bgc
\caption{Results from the coefficient matrix 'atmosmodj'.}
\edc
\includegraphics[scale=0.3]{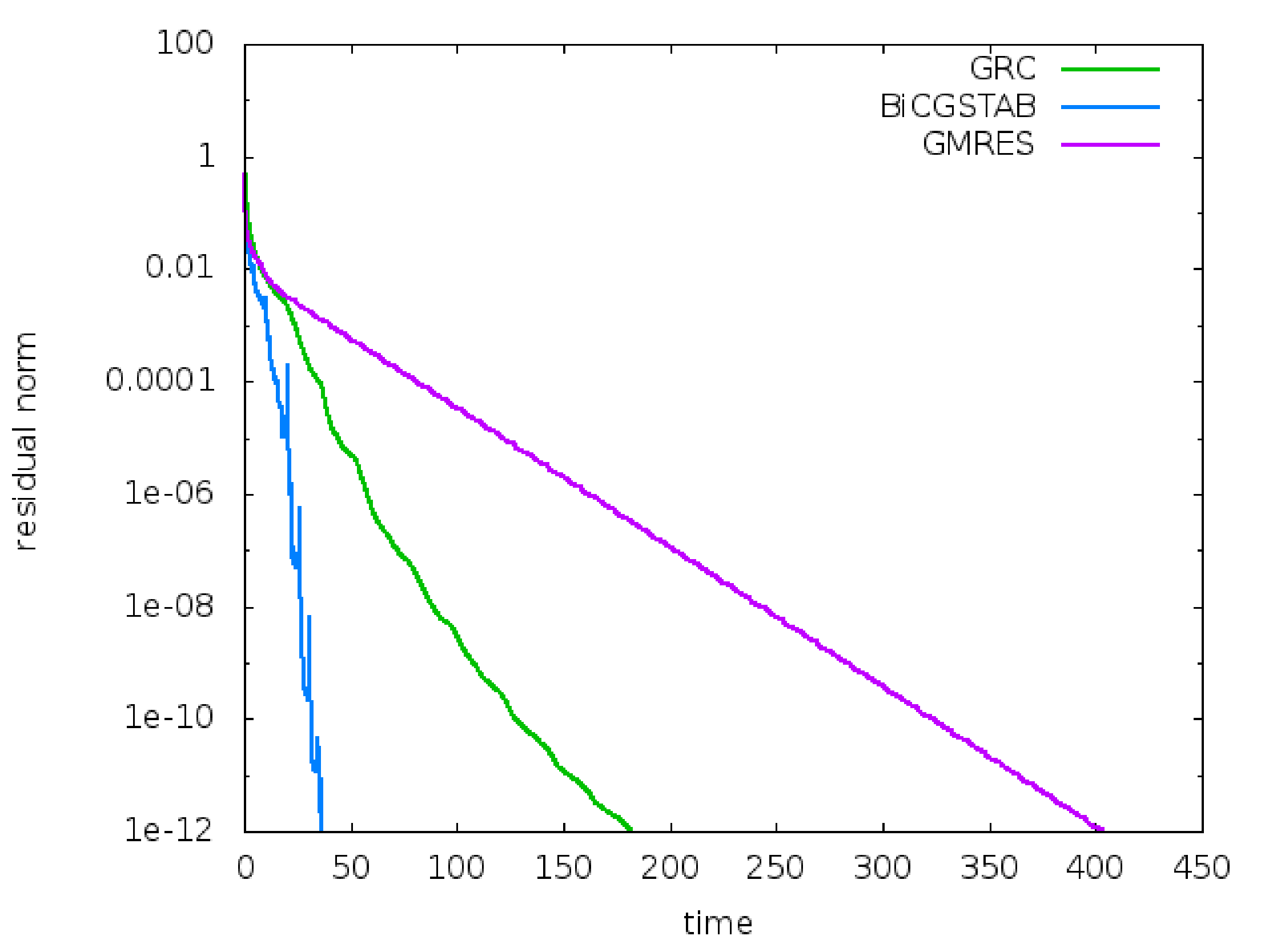}
\label{fig3.8}

\end{figure}
\begin{figure}[!h]
  
\bgc

\vspace{3mm}

\caption{Results from the coefficient matrix 'ML\_Laplace'.}
\edc
\includegraphics[scale=0.3]{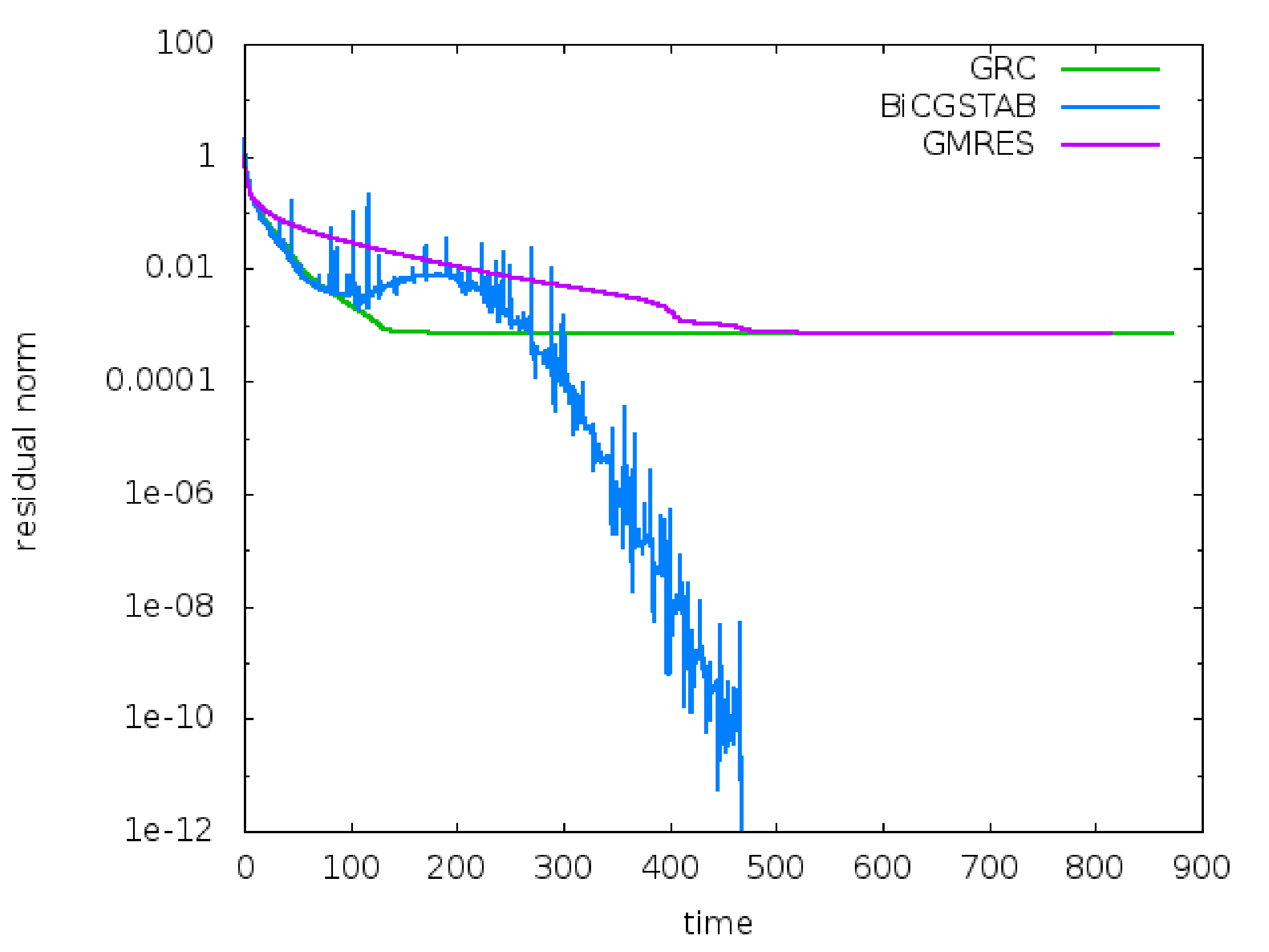}
\label{fig3.9}

\end{figure}
\begin{figure}[!h]

\bgc
  
\vspace{3mm}

\caption{Results from the coefficient matrix 'Transport'.}
\edc
\includegraphics[scale=0.3]{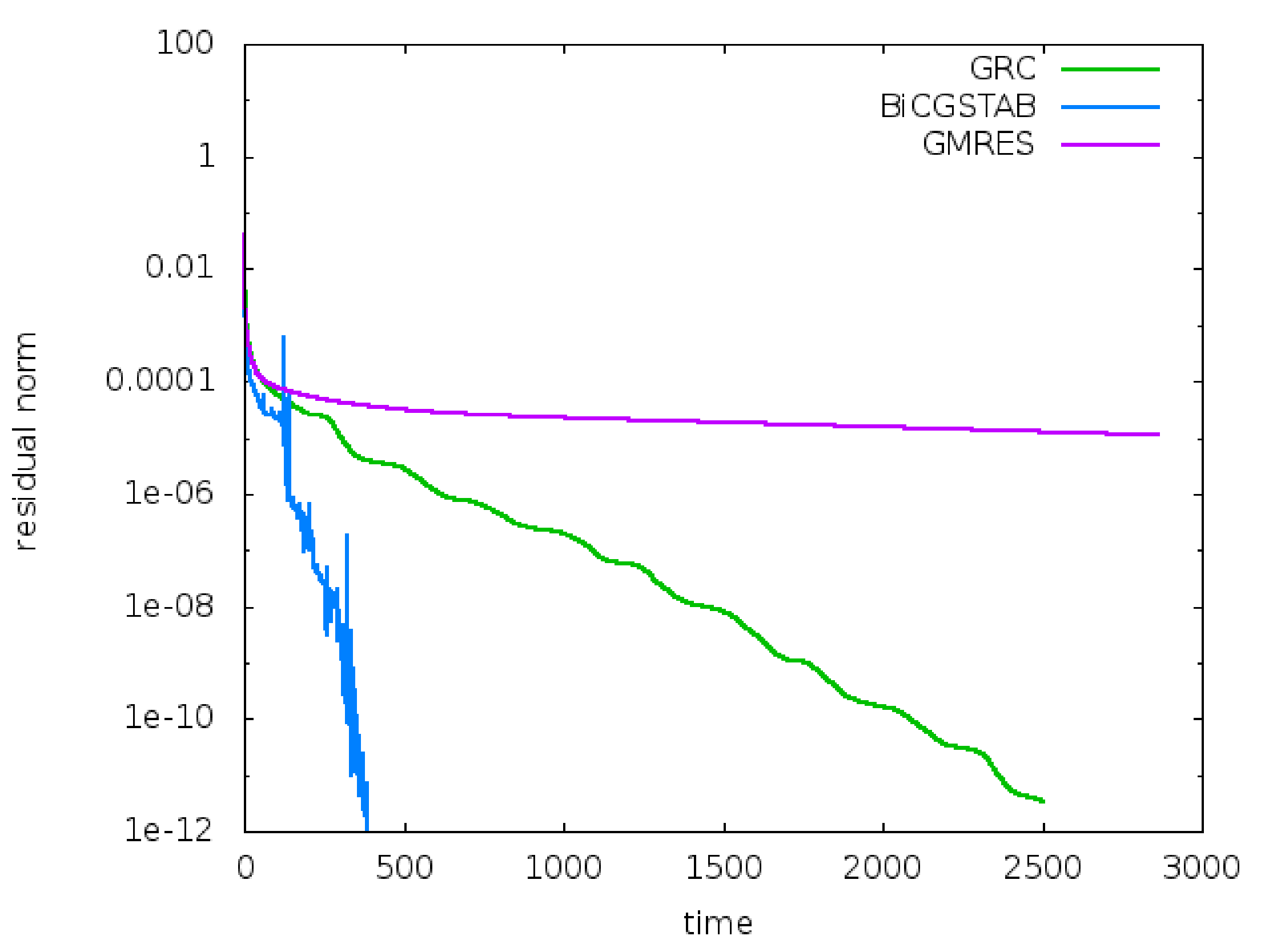}
\label{fig3.10}
\end{figure}

%\setlength\fboxsep{0pt}
%\setlength\fboxrule{0.5pt}
%\fbox{\includegraphics{iter_standard.png}}

We have shown that GRC is a Krylov subspace method and
its close relationship to the conjugate residual method.
Also,  numerical experiments indicate that it works as well as GMRES and BiCGSTAB
for general unsymmetric sparse problems.
Among the four methods reported in this paper, the RC method is
considered to be most effective with matrices
for which a relaxation method works effectively.
Among the Krylov subspace methods, no single method has been shown to be superior to others,
however,
the GRC and GMRES methods, show robust convergence by the residual norm minimization.
In addition, GRC $(L=5)$, needs to keep less number of vectors than
GMRES $(K=40)$, thus we expect that GRC has an advantage for
significantly larger matrix sizes.

\end{document}